\documentclass[11pt,legal]{article}
\pdfoutput=1
\usepackage{theorem,ifthen}
\usepackage{algorithm}
\usepackage{algorithmic}
\usepackage{amssymb,amsfonts,amsmath,latexsym,dsfont}
\usepackage{tikz,pgflibraryplotmarks}
\usepackage{fullpage}
\usepackage{graphicx}
\usepackage{mathrsfs}
\usepackage{authblk}
\usepackage{gnuplot-lua-tikz}
\usepackage{url}
\usepackage{booktabs}  % Nice tabulars

\graphicspath{{images/}{./}}

\usepackage{boxedminipage}

\usepackage{pgf,pgfarrows} %% pdf-drawing package
\usepackage{subcaption} 

\newcommand{\FrameboxA}[2][]{#2}
\newcommand{\Framebox}[1][]{\FrameboxA}

\newcommand{\bfA}{{\bf A}}
\newcommand{\bfB}{{\bf B}}

\newcommand{\bfE}{{\bf E}}

\newcommand{\bfG}{{\bf G}}

\newcommand{\bfI}{{\bf I}}

\newcommand{\bfK}{{\bf K}}
\newcommand{\bfL}{{\bf L}}

\newcommand{\bfR}{{\bf R}}

\newcommand{\bfU}{{\bf U}}
\newcommand{\bfV}{{\bf V}}
\newcommand{\bfW}{{\bf W}}

\newcommand{\bfPhi}{{\bf \Phi}}

\newcommand{\bfmu}{{\boldsymbol{\mu}}}
\newcommand{\bftheta}{{\boldsymbol{\theta}}}

\newcommand{\bfc}{{\bf c}}
\newcommand{\bfe}{{\bf e}}

\newcommand{\bfy}{ {\bf y}}
\newcommand{\bfu}{{\bf u}}

\newcommand{\bfz}{{\bf z}}
\newcommand{\bfzero}{{\boldsymbol{0}}}

\newcommand{\R}{\ensuremath{\mathds{R}}}

\newlength\myindent
\theoremstyle{remark}
\newtheorem{remark}{Remark}

\title{Layer-Parallel Training of Deep Residual Neural Networks}

\author[1]{S. G\"unther\thanks{Corresponding author: stefanie.guenther@scicomp.uni-kl.de}}
\author[2]{L. Ruthotto\thanks{LR is supported by the US National Science Foundation awards DMS 1522599 and DMS 1751636}}
\author[3]{J.B. Schroder}
\author[4]{E.C. Cyr\thanks{Sandia National Laboratories is a multimission laboratory managed and operated by
National Technology \& Engineering Solutions of Sandia, LLC, a wholly owned
subsidiary of Honeywell International Inc., for the U.S. Department of Energy's
National Nuclear Security Administration under contract DE-NA0003525.
The views expressed in the article do not necessarily represent the views of the U.S.
Department of Energy or the United States Government.}}
\author[1]{N.R. Gauger}

\affil[1]{Scientific Computing Group, TU Kaiserslautern, Germany}
\affil[2]{Department of Mathematics and Computer Science, Emory University, Atlanta, GA, USA}
\affil[3]{Dept. of Mathematics and Statistics, University of New Mexico, Albuquerque, NM, USA}
\affil[4]{Computational Mathematics Department, Sandia National Laboratories, Albuquerque, NM, USA}

\begin{document}

\maketitle

\begin{abstract}
Residual neural networks (ResNets) are a promising class of deep neural networks that have shown excellent performance for a number of learning tasks, e.g., image classification and recognition.
Mathematically, ResNet architectures can be interpreted as forward Euler discretizations of a nonlinear initial value problem whose time-dependent control variables represent the weights of the neural network. 
Hence, training a ResNet can be cast as an optimal control problem of the associated dynamical system. 
For similar time-dependent optimal control problems arising in engineering applications, parallel-in-time methods have shown notable improvements in scalability.
This paper demonstrates the use of those techniques for efficient and effective training of ResNets. 
The proposed algorithms replace the classical (sequential) forward and backward propagation through the network layers by a parallel nonlinear multigrid iteration applied to the layer domain. 
This adds a new dimension of parallelism across layers that is attractive when training very deep networks.
From this basic idea, we derive multiple layer-parallel methods. 
The most efficient version employs a simultaneous optimization approach where updates to the network parameters are based on inexact gradient information in order to speed up the training process.
Using numerical examples from supervised classification, we demonstrate  that the new approach achieves similar training performance to traditional methods, but enables layer-parallelism and thus provides speedup over layer-serial methods through greater concurrency. 
\end{abstract}

%%%%%%%%%%%%%%%%%%%%%%%%%%%%%%%%%%%%%%%%%%%%%%%%%%%%%%%%%%%%%%%%%%%%%%
\section{Introduction}
\label{sec:introduction}
%%%%%%%%%%%%%%%%%%%%%%%%%%%%%%%%%%%%%%%%%%%%%%%%%%%%%%%%%%%%%%%%%%%%%%
One of the most promising areas in artificial intelligence is deep learning, a form of machine learning that uses neural networks containing many hidden layers~\cite{bengio2009learning,lecun2015deep}.  
Deep neural networks (DNNs), and in particular deep residual networks (ResNets)~\cite{he2016deep}, have been breaking human records in various contests and are now central to technology such as image recognition~\cite{hinton2012deep,KrizhevskySutskeverHinton2012,lecun2015deep} and natural language processing~\cite{BordesEtAl2014,CollobertEtAl2011,  JeanEtAl2014}. 

The abstract goal of machine learning is to model a function $f : \R^n \times \R^p \to \R^m$ and train its parameter $\bftheta\in\R^p$ such that 
\begin{equation}\label{eq:interp}
f(\bfy, \bftheta) \approx \bfc
\end{equation}
for input-output pairs $(\bfy,\bfc)$ from a certain data set $\mathcal{Y} \times \mathcal{C}$. 
Depending on the nature of inputs and outputs, the task can be regression or classification. When outputs are available for all samples, parts of the samples, or are not available, this formulation describes supervised, semi-supervised, and unsupervised learning, respectively. 
The function $f$ can be thought of as an interpolation or approximation function.

In deep learning, the function $f$ involves a DNN that aims at transforming the input data using many layers. 
The layers successively apply affine transformations and element-wise nonlinearities that are parametrized by the network parameters $\bftheta$. 
The training problem consists of finding the parameters $\bftheta$ such that~\eqref{eq:interp} is satisfied for data elements from a training data set, but also holds for previously unseen data from a validation data set, which has not been used during training. 
The former objective is commonly modeled as an expected loss and optimization techniques are used to find the parameters that minimize the loss.

Despite rapid methodological developments, compute times for training state-of-the-art DNNs can still be prohibitive, measured in the order of hours or days, involving hundreds or even thousands of layers and millions or billions of network parameters~\cite{dean2012large, keuper2016distributed}. 
There is thus a great interest in increasing parallelism to reduce training runtimes.
The most common approach involves \textit{data-parallelism}, where elements of the training data set are distributed onto multiple compute units. Synchronous and asynchronous data-parallel training algorithms have been developed to coordinate the network parameter updates~\cite{iandola2016firecaffe,abadi2016tensorflow}. 
Another approach is referred to as \textit{model-parallelism}, which aims at partitioning different layers of the network and its parameters to different compute units. Model parallelism has traditionally been used when the network dimension exceeds available memory of a single compute unit. Often, a combination of both approaches is employed~\cite{harlap2018pipedream, dean2012large}. 
% When analyzed, these approaches solve a modified
% optimization problem when compared to a straight-forward serial processing
% of all training batches, with the result that training effectiveness can suffer.  In contrast,
% the layer parallelism considered here does not change the serial processing of
% training batches.

However, none of the above approaches to parallelism tackle the scalability barrier created by the intrinsically serial propagation of data through the network itself. In either of the above approaches, each subsequent layer can process accurate information only after the previous layer has finished its computation. As a result, training runtimes typically scale linearly with the number of layers. 
As current state-of-the-art networks tend to increase complexity by adding more and more layers (see, e.g., the ResNet-1001 with 1001 layers and 10.2 million weights in~\cite{He2016identity}), the serial layer propagation creates a serious bottleneck for fast and scalable training algorithms seeking to leverage modern HPC facilities. 

In this paper, we address the above scalability barrier by introducing concurrency across the network layers. To this end, we replace the serial data propagation through the network layers by a nonlinear multigrid method that treats layers, or layer chunks, simultaneously and thus enables full layer-parallelism. 
Our goal is to have a training methodology that is scalable in the number of layers, e.g., doubling the number of layers and the number of compute resources should result in a near constant runtime.
To achieve this, we leverage recent advances in parallel-in-time integration methods for unsteady differential equations.\footnote{For an introduction and overview on various parallel-in-time integration schemes for unsteady differential equations we refer the reader to the review paper in \cite{Ga2015}, and more recent development such as \cite{GanderGuttelPetcu2018, gotschel2019efficient}}

The forward propagation through a ResNet can be seen as a discretization of a time-dependent ordinary differential equation (ODE), which was first observed in~\cite{HaberRuthotto2017,E2017,HaberHolthamRuthotto2017}. Interpreting the network propagation as a nonlinear dynamical system has since attracted increasing attention (see, e.g. \cite{lu2017beyond} and references therein, or \cite{chen2018neural}).

Based on this interpretation, we employ a multigrid reduction in time approach \cite{Fa2014} that divides the time domain -- which, in this interpretation, corresponds to the layer domain -- into multiple time chunks that can be processed in parallel on different compute units. Coupling of the chunks is achieved through a coarse-grid correction scheme that propagates information across chunk interfaces on a coarser time- (i.e. layer-) grid.   
The method can be interpreted as a parallelization of the model, processing layer chunks simultaneously within the iterative multigrid scheme thus breaking the traditional layer-serial propagation. At convergence, the iterative multigrid scheme solves the same problem as a layer-serial method and it can thus be utilized in any common gradient-based optimization technique to update the network parameters, such as the stochastic gradient descent (SGD) or other batch approaches, without loss of accuracy. Further, it can be applied in addition to any data-parallelism across the data set elements, thus multiplying data-parallel runtime speedup. Runtime speedup over traditional layer-serial methods is achieved through the new dimension of parallelization across layers enabling greater concurrency.

The addition of layer-parallelism
allows for ResNets to take advantage of large machines currently
programmed with message-passing style parallelism.  The use of such large
machines in conjunction with a multilevel training algorithm scalable in
the number of layers, opens the door to training networks with thousands,
or possibly even millions of layers.
We demonstrate the feasibility of such an approach by
using the parallel-in-time package XBraid~\cite{xbraid-package} with a ResNet on large
clusters.  Additionally, the non-intrusive approach of XBraid would allow
for any node-level optimizations (such as those utilizing GPUs) to be
used. 

The iterative nature of the multigrid approach further enables the use of simultaneous optimization algorithms for training the network. Simultaneous optimization methods have been widely used for optimization problems that are constrained by partial differential equations (PDEs), where they show promise for reducing the runtime overhead of the optimization when compared to a pure simulation of the underlying PDE (see e.g. ~\cite{borzi2011computational, ziems2011adaptive, biros2005paralleli} and references therein). They aim at solving the optimization problem in an all-at-once fashion, updating the optimization parameters simultaneously while solving for the time-dependent system state. 
Here, we apply the \textit{One-shot} method ~\cite{bosse2014oneshot, gunther2018xbraid} to solve the training problem simultaneously for the network state and parameters. In this approach, network parameter updates are based on inexact gradient information resulting from early stopping of the layer-parallel multigrid iteration.

The paper is structured as follows: Section \ref{sec:oc} gives an introduction to the deep learning optimization problem and its interpretation as an optimal control problem. Further, it discusses numerical discretization of the optimal control problem and summarizes necessary conditions for optimality. 
We then introduce the layer-parallel multigrid approach replacing the forward and backward propagation through the network in Section \ref{sec:pinl}.
Section \ref{sec:oneshot} focuses on the integration of the layer-parallel multigrid scheme into a simultaneous optimization algorithm. 
Numerical results demonstrating the feasibility and runtime benefits of the proposed layer-parallel scheme are presented in Section \ref{sec:numerics}.
 
%%%%%%%%%%%%%%%%%%%%%%%%%%%%%%%%%%%%%%%%%%%%%%%%%%%%%%%%%%%%%%%%%%%%%%
\section{Deep Learning as a Dynamic Optimal Control Problem}
\label{sec:oc}
%%%%%%%%%%%%%%%%%%%%%%%%%%%%%%%%%%%%%%%%%%%%%%%%%%%%%%%%%%%%%%%%%%%%%%

In this section, we present an optimal control formulation of a supervised classification problem using a deep residual network.
Limiting the discussion to this specific task allows us to provide a self-contained mathematical description.
We note that the layer-parallel approach can be extended to other learning tasks, e.g., semi-supervised learning, auto-regression, or recurrent learning and refer to~\cite{GoodfellowEtAl2016, abu2012learning} for a general introduction and a comprehensive overview of deep learning techniques.

\subsection{Optimal Control Formulation}
In supervised classification, the given data set consists of $s$ feature, or example, vectors $\bfy_1,\bfy_2,\ldots, \bfy_s \in \R^{n_f}$ and associated class probability vectors $\bfc_1, \bfc_2, \ldots, \bfc_s \in \Delta_{n_c}$, where $\Delta_{n_c}$ denotes the unit simplex in $\R^{n_c}$, and $n_f, n_c \in \mathds{N}$ denote the number of features and classes in the given data set, respectively. 
The $j$-th component of $\bfc_k$ represents the probability of example $\bfy_k$ belonging to the $j$-th class.
The learning problem aims at training a network, and its classifier, that approximate the feature-to-class mapping for all data elements.

A powerful class of networks are Residual Neural Networks (ResNets)~\cite{he2016deep}. 
In an abstract form, the network transformation to a generic input data example $\bfy$ using an $N$-layer ResNet can be written as
\begin{equation}\label{eq:resnetpropagation}
	\bfu^{n+1} = \bfu^n + h F(\bfu^n, \bftheta^n), \quad \text{ for } n=0,1,\ldots,N-1, \quad \text{with} \quad \bfu_0 = \bfL_{\rm in}\bfy.
\end{equation}
with $\bfu^n\in\R^q$, $q$ being the network width.
The transformations in $F$ typically consist of affine linear and element-wise nonlinear transformations that are parameterized by the entries in the layer weights $\bftheta^{0},\ldots, \bftheta^{N-1} \in \R^d$, respectively. For simplicity, we consider the single layer perceptron model
\begin{equation}
	\label{eq:singlelayerperceptron}
	F(\bfu,\bftheta) = \sigma( \bfK(\bftheta^{(1)}) \bfu + \bfB \bftheta^{(2)}), 
\end{equation}
where $\sigma : \R \to \R$ is a nonlinear activation function that is applied component-wise, e.g., $\sigma(x) = \tanh(x)$ or $\sigma(x) = \max\{x,0\}$.
Here, each weight vector $\bftheta \in \{ \bftheta^0, \ldots, \bftheta^{N-1}\}$ is partitioned into one part that defines a linear operator $\bfK\in \R^{q\times q}$ and another part that represents coefficients of a bias with respect to columns of the given matrix $\bfB \in \R^{q\times n_b}$ (e.g., $n_b=1$ and $\bfB = \bfe_q$, a vector of all ones, to add a constant shift to all features).
In this work, we assume that the linear operators $\bfK(\cdot)$ are either dense matrices or correspond to convolutional operators (see~\cite{LeCunEtAl1990}) parametrized by $\bftheta^{(1)}$, whose entries we determine in the training. However, our method can be extended to other parameterizations (e.g., the layer used in~\cite{he2016deep}, which features two affine transformations, or the layer based on an antisymmetric matrix suggested in ~\cite{HaberRuthotto2017}).
While $\bfK(\cdot)$ needs to be a square matrix we use a non-square model for the operator $\bfL_{\rm in} \in \R^{q\times n_f}$ to map the data set elements to the network width.

Considering a small but positive $h$ in \eqref{eq:resnetpropagation}, it is intuitive to interpret the ResNets propagation \eqref{eq:resnetpropagation} as a forward Euler discretization of the initial value problem
\begin{equation}\label{eq:ResNNcont}
	\partial_t \bfu(t) = F(\bfu(t), \bftheta(t)), \quad t \in [0,T], \quad \text{with} \quad  \bfu(0) = \bfL_{\rm in}\bfy.
\end{equation}
In this formulation, $t$ is an artificial time that refers to the propagation of the input features through the neural network, starting from the input layer with $\bfu(0)$ to the network output $\bfu(T)$ being the solution of the initial value problem evaluated at some final time $T$. In contrast to the discrete ResNet propagation, the dynamical system continuously transforms the network state $\bfu(t)$ by prescribing its time-derivative with the vector field $F$, whose parameters $\bftheta(t)$ will be learned during training.

In order to classify the network output into a specific class, a hypothesis function is required that predicts the class probabilities. 
Here, we limit ourselves to multinomial regression models, which are common in deep learning. To this end, we consider the \textit{softmax hypothesis function} given by
\begin{equation}\label{eq:softmax}
	{S}(\bfu(T),\bfW,\bfmu) = \frac{1}{\bfe_{n_c}^\top \exp(\bfz)} {\exp(\bfz)}, \quad \bfz = \bfW \bfu(T) + \bfmu,
\end{equation}
where $\bfe_{n_c} \in \R^{n_c}$ is a vector of all ones, $\exp$ is the exponential function applied element-wise, and $\bfW \in \R^{n_c \times q}$, $\bfmu \in \R^{n_c}$ denote a weight matrix and bias vector, whose entries need to be learned in training alongside the network parameters $\bftheta(t)$. 

The performance of the network transformation and classification can then be measured by comparing the predicted class probabilities to the given ones in $\bfc$. To this end, we use the \textit{cross entropy loss function}
\begin{align}\label{eq:loss}
	\ell(\bfu(T), \bfc, \bfW,\bfmu) & =  -  \bfc^\top \log(S(\bfu(T),\bfW,\bfmu)).
	% \\ &=  -  \bfc^\top  \bfz + \log(\bfe_{n_c}^\top \exp(\bfz)).
\end{align}

The training problem consists of minimizing the average cross entropy loss function over many examples with respect to $\bftheta(t)$,$\bfW$ and $\bfmu$. It can thus be cast as the following continuous-in-time optimal control problem:
\begin{align}\label{eq:opt}
	\min %_{\bfu_1,\ldots,\bfu_s,\bftheta,\bfW,\bfmu}  
	\quad  \frac1s \sum_{k=1}^s \ell(&\bfu_k(T), \bfc_k, \bfW,\bfmu) +  \int_0^T R(\bftheta(t), \bfW, \bfmu) \, \mathrm{d}t \\
	\text{ subject to }\quad  \partial_t \bfu_k(t) &= F(\bfu_k(t), \bftheta(t)), \quad \forall t \in [0,T], \label{eq:dyn} \\
	\bfu_k(0) &= \bfL_{\rm in}\bfy_k, \hspace{11ex} \forall \, k = 1, \dots, s.
	\label{eq:dyn_init}
\end{align}
The optimal control problem aims at finding control variables $\bftheta(t)$ (the network weights) and corresponding state variables $\bfu_k(t)$ (the network states) that minimize the objective function, while satisfying the constraints \eqref{eq:dyn}--\eqref{eq:dyn_init} (the network dynamics) for all data set elements $k=1,\ldots, s$. The objective function consists of the empirical cross entropy loss function evaluated at the final time $T$, and an additional regularization term denoted by $R$.
In conventional deep learning approaches, $R$ typically applies a Tikhonov regularization penalizing ``large'' network and classification parameters, measured in a chosen norm. Within the time-continuous optimal control interpretation, we additionally penalize the time-derivative of $\bftheta(t)$ in order to ensure weights that vary smoothly in time. This is an important ingredient for stability analysis~\cite{HaberRuthotto2017}. 

\subsection{Discretization of the Optimal Control Problem}

We solve the continuous-in-time optimal control problem in a first-discretize-then-optimize fashion. 
We discretize the control $\bftheta(t)$ and the states $\bfu(t)$ at regularly spaced time points $t_n = n\cdot h$, where $h = T/N$ and $n=0,1,\ldots,N$. In this setting, each discrete state $\bfu^n$ and control $\bftheta^n$ correspond to the $n$-th layer of the network.
 This leads to the discrete control problem
\begin{align}\label{eq:disc_opt}
	\min \quad  \frac1s \sum_{k=1}^s \ell(&\bfu_k^N, \bfc_k, \bfW,\bfmu) + \sum_{n=0}^{N-1} R(\bftheta^n, \bfW, \bfmu) \\
   \label{eq:disc_opt:forward_problem}
	\text{ subject to }  \quad  \bfu_k^{n+1} &= \bfPhi(\bfu_k^n,\bftheta^n), \quad \forall \, n = 0,\dots, N-1 \\ 
   \label{eq:disc_opt:forward_problem_init}
	\bfu_k^0 &= \bfL_{\rm in}\bfy_k \hspace{8ex} \forall \, k \in 1, \dots, s.
\end{align}
In this general description, $\bfPhi$ can denote any layer-to-layer propagator which maps data $\bfu^n$ to the next layer. In case of a forward Euler discretization, it reads
\begin{equation}\label{eq:Phi}
	\bfPhi(\bfu^n,\bftheta^n) = \bfu^n + h F(\bfu^n,\bftheta^n),
\end{equation}
giving the ResNet propagation as in \eqref{eq:resnetpropagation}. However, the time-continuous interpretation of the network propagation permits to employ other, possibly more stable discretization schemes of the initial value problem \eqref{eq:ResNNcont}, see \cite{HaberRuthotto2017}, and thus opens the door for new network architecture designs.
It also allows for discretization of the controls and states at different time points, which can improve the efficiency and is the subject of further research. 
Further, numerical advances for solving the corresponding optimal control problem can be leveraged, such as the time-parallel approach which is discussed in this paper.

\subsection{Necessary Optimality Conditions}

The necessary conditions for optimality of the discrete, equality-constrained optimization problem \eqref{eq:disc_opt}-\eqref{eq:disc_opt:forward_problem_init} can be derived from the associated Lagrangian function
\begin{align}\label{eq:Lagrangian}
	L :=  &J  
		+\sum_{k=1}^s \left[ \sum_{n=0}^{N-1}  \left(\bar \bfu_k^{n+1}\right)^T\left( \bfPhi(\bfu^n_k, \bftheta^n) - \bfu^{n+1}_k \right) + \left(\bar \bfu^0_k\right)^T\left(\bfL_{\rm in}\bfy_k - \bfu^0_k\right) \right]
\end{align}
where $J$ denotes the objective function in \eqref{eq:disc_opt} and $\bar \bfu^n_k$ are called adjoint variables for layer $n=0,\dots,N$ and example $k=1,\dots,s$. Optimal points of the problem are saddle points of the Lagrangian function (see e.g. \cite{nocedal2006numerical}), thus equating its partial derivatives with respect to all state, adjoint and control variables  
% $\bfu^0,\dots, \bfu^N$, $\bar \bfu^0 \dots, \bar \bfu^N$, and $\bftheta, \bfW, \bfmu$
to zero yields the following necessary conditions for optimality:   
\begin{enumerate}
	\item \textit{State equations}
		\begin{align}
			\bfu_k^{n+1} &=\bfPhi(\bfu_k^n,\bftheta^n), \quad \forall n = 0,\dots, N-1, \label{eq:discretestate1}\\
			  \text{with} \quad  \bfu_k^0 &= \bfL_{\rm in}\bfy_k, \qquad \qquad	\text{for all } k=1,\dots,s \label{eq:discretestate2}
		\end{align}
	\item \textit{Adjoint equations}
		\begin{align}
			\bar \bfu^{n}_k &= \left(\partial_{\bfu} \bfPhi(\bfu^n_k,\bftheta^n)\right)^T \bar \bfu^{n+1}_k ,\quad \forall n=0, \ldots,N-1, \label{eq:discreteadjoint1}\\
			 \text{with} \quad \bar \bfu^{N}_k &= \frac 1s \left(\partial_{\bfu} \ell(\bfu^N_k,\bfc_k, \bfW, \bfmu)\right)^T  \quad	\text{for all } k=1,\dots,s \label{eq:discreteadjoint2}
		\end{align}
	\item \textit{Design equations}
		\begin{align}
		 0 &=  \sum_{k=1}^s\left(\partial_{\bftheta^n}\bfPhi(\bfu^n_k,\bftheta^n)\right)^T\bar \bfu_k^{n+1} + \left(\partial_{\bftheta^n} R\right)^T \quad	\forall n=0,\dots,N-1 \label{eq:reducedgradient1} \\
		 0 &=	 \frac 1s \sum_{k=1}^s \left(\partial_{ \bfW,\bfmu}\ell(\bfu^N_k, \bfc_k, \bfW, \bfmu)\right)^T + \left(\partial_{\bfW,\bfmu} R\right)^T  \label{eq:reducedgradient2}
		\end{align}
\end{enumerate}
Here, subscripts denote partial derivatives, $\partial_x = \frac{\partial}{\partial x}$. 
Training a residual network corresponds to the attempt of solving the above set of equations for the special choice of $\bfPhi$ being the forward Euler time-integration scheme. However, the above equations, as well as the discussions in the remainder of this paper, are general, in the sense that any other layer-to-layer propagator $\bfPhi$ can be utilized that corresponds to the discretization of the dynamical system \eqref{eq:ResNNcont}.

The state equations correspond to the forward propagation of input examples $\bfy_k$ through the network layers. The adjoint equations propagate partial derivatives with respect to the network states backwards through the network layers, starting from a terminal condition at $N$ equal to the local derivative of the loss function. 
In a time-continuous setting, the adjoint equations correspond to the discretization of an additional adjoint dynamical system for propagating network state derivatives backwards in time.\footnote{The adjoint approach is a common and well-established method in optimal control that provides gradient information at computational costs that are independent of the design space dimension, see e.g. \cite{Giles2000}}
We note that solving the above adjoint equations backwards in time is equivalent to the \textit{backpropagation} method that is established within the deep learning community for computing the network gradient~\cite{LeCunEtAl1990}. It further corresponds to the reverse mode of automatic differentiation~\cite{griewank2008evaluating}.
The adjoint variables are utilized in the right-hand-side of the design equations, which then form the so-called reduced gradient. For feasible state and adjoint variables, the reduced gradient holds the total derivative, i.e. the sensitivity, of the objective function with respect to the controls. It is thus used within gradient-based optimization methods for updating the network controls.

%%%%%%%%%%%%%%%%%%%%%%%%%%%%%%%%%%%%%%%%%%%%%%%%%%%%%%%%%%%%%%%%%%%%%%
\section{Layer-Parallel Multigrid Approach}
\label{sec:pinl}
%%%%%%%%%%%%%%%%%%%%%%%%%%%%%%%%%%%%%%%%%%%%%%%%%%%%%%%%%%%%%%%%%%%%%%

In order to achieve concurrency across all the network layers, we
replace the sequential propagation through the residual network (forward and backward)
with an iterative multigrid scheme. 

Based on the time-continuous nonlinear ODE interpretation of ResNets as in \eqref{eq:ResNNcont}, and its time-discretization as in \eqref{eq:disc_opt:forward_problem}--\eqref{eq:disc_opt:forward_problem_init}, we employ the multigrid reduction in time (MGRIT)~\cite{Fa2014} method to parallelize across the time domain of the network. 
While the discussion in this section revolves around time-grids, here each time point
is considered a layer in the network.  Thus, the multigrid approach constructs a multilevel
hierarchy, where each level is a network containing fewer layers (i.e., fewer time points).  The
coarsest level will contain only a handful of layers, while the finest level
could contain thousands (or more) of layers.
When run in parallel, each compute unit will own only a few {fine-grid} layers,
thus allowing for massive parallelism to be applied to the learning algorithm.

The MGRIT scheme was introduced in
\cite{Fa2014} and first applied to neural networks in \cite{Sc2017}, although
that work considered parallelism over epochs of the training algorithm, not layers.  
We refer to the works
\cite{Fa2014, falgout2017multigrid, guenther2017xbraid} for the details of the 
method, but we will here provide a self-contained overview of the MGRIT scheme.

\subsection{Multigrid Across Layers for Forward Propagation}
\label{sec:mgrit}
Consider the network states to be collected in a vector $\bfU = (\bfu^0, \bfu^1, \ldots, \bfu^N)$. The forward propagation through the network \eqref{eq:discretestate1}--\eqref{eq:discretestate2} can than be written as the space-time system
\begin{align}\label{eq:nonlinear_mgrit}
	\bfA(\bfU,\bftheta) := 
	\begin{pmatrix} 
		\bfu^0 \\
		\bfu^1 - \bfPhi(\bfu^0, \bftheta^0) \\
		\vdots \\
		\bfu^{N} - \bfPhi(\bfu^{N-1}, \bftheta^{N-1}) \\
	\end{pmatrix} 
	=
	\begin{pmatrix}
		\bfL_{\rm in}\bfy\\
		\bfzero \\
		\vdots \\
		\bfzero 
	\end{pmatrix}
	=: \bfG.
\end{align}
where each block row corresponds to a time step, which in turn corresponds to a layer in the network. 
Here, the $\bfu^n$ denote the network states at each time step for either a single generic input vector $\bfy$ or for a batch, i.e. a subset, of input vectors $\bfy_k,\, k \in \mathcal S\subset \{1,\dots,s\}$.

Sequential time
stepping solves (\ref{eq:nonlinear_mgrit}) through forward substitution, i.e. forward propagation of input data through the network layers.
In contrast, MGRIT solves (\ref{eq:nonlinear_mgrit}) iteratively, beginning with some initial solution guess for $\bfU$, by using  
the Full Approximation Storage (FAS) nonlinear multigrid method \cite{Brandt_1977}, see Section \ref{sec:FAS}. In both cases, the exact same equations are solved and thus the same solution is reached 
(in the case of MGRIT, to within a user tolerance). Regarding cost, sequential time-stepping is $O(N)$, but sequential.  Instead, MGRIT solves this system with an $O(N)$ multigrid method with a larger
computational constant, but with parallelism in the layer dimension.  This parallelism allows for a
distributed workload, processing multiple layers in parallel on multiple compute units.
Typically, a certain number of processors are needed for MGRIT to show
a speedup over layer-serial forward propagation.  This is referred to as the cross-over point. 
However, the speedups 
observed can be large, e.g., the work \cite{gunther2018xbraid} showed a speedup
of 19x for a model optimization problem while using an additional 256 processors in time.

\subsubsection{MGRIT using Full Approximation Scheme (FAS)}
\label{sec:FAS}

\begin{figure}
\center
\usetikzlibrary{decorations.pathreplacing}
\begin{tikzpicture}
\path[draw] (0,0) -- (10,0);
%%
% Draw Fine Points
\foreach \x in {0,1,...,3}
   \path[draw, line width = 1pt] (\x/2,0.1) -- (\x/2,-0.1) node [below] {$t_{\x}$};
\path[draw, line width = 1pt] (2,0.1) -- (2,-0.1) node [below] {$...$};
\path[draw, line width = 1pt] (2.5,0.1) -- (2.5,-0.1) node [below] {$t_c$};
\foreach \x in {4,5,...,20}
   \path[draw, line width = 1pt] (\x/2,0.1) -- (\x/2,-0.1);
%%
% Draw Coarse Points
\filldraw[color = red] (0,0) circle (3.5pt);
   \filldraw[color = red] (0.05,0.1) node [above] {$T_0$};
\filldraw[color = red] (2.5,0) circle (3.5pt);
   \filldraw[color = red] (2.95,0.1) node [above] {$T_1 \;\; \dots $};
\filldraw[color = red] (5,0) circle (3.5pt);
\filldraw[color = red] (7.5,0) circle (3.5pt);
\filldraw[color = red] (10,0) circle (3.5pt);
   \filldraw[color = red] (10.25,0.1) node [above] {$T_{N/c}$};
%%
% Draw Legend
\filldraw[color = red] (11.7,.5) circle (3.5pt);
   \filldraw[color = black] (12.7,0.5) node {: C-point};
   \filldraw[color = black] (13.65,0.14) node {\footnotesize{(fine and coarse grid)}};
\path[draw, line width = 1pt] (11.4,-0.5) -- (11.8,-0.5);
\path[draw, line width = 1pt] (11.6,-0.37) -- (11.6,-0.63);
   \filldraw[color = black] (12.7,-0.5) node {: F-point};
   \filldraw[color = black] (13.25,-0.86) node { \footnotesize{(fine grid only)}};
   \draw [decoration={brace,mirror}, decorate, line width = 1pt] (5,-0.8) -- (7.5,-0.8) node [below, pos=0.5] {$h_{\Delta} = c h$};
\draw [decoration={brace,mirror}, decorate, line width = 1pt] (5.5,-0.22) -- (6,-0.22) node [below, pos=0.5] {$h $};
\end{tikzpicture}
\caption{Fine grid ($t_i$) and coarse grid ($T_j$) for coarsening factor
   $c=5$.  MGRIT eliminates the fine points (black vertical lines, F-points) 
   to yield a coarse level composed of the red circles (C-points).} 
\label{fig-grid}
\end{figure}
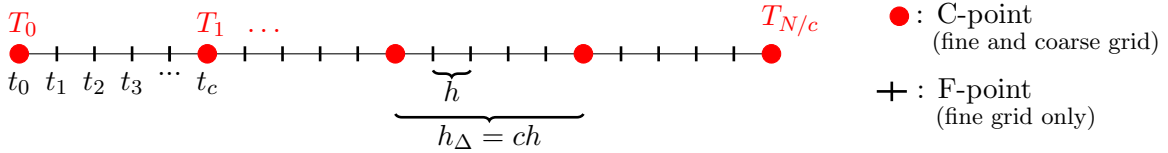

Similar to linear multigrid methods, the nonlinear FAS method computes coarse-grid error corrections to fine-grid approximations of the solution.
Each iteration of the nonlinear MGRIT scheme consists of three steps: First, a relaxation scheme is employed to cheaply compute an approximation to the true solution on the fine grid. Then, the current error is approximated on a coarser grid by solving a coarse-grid residual equation. Lastly, the interpolated coarse-grid error approximation is used to correct the current fine-grid solution approximation. 
This idea is based on the fact that low frequency error components can be reduced with relaxation much faster on coarser grids. While a general introduction to linear and nonlinear multigrid methods can be found in \cite{BrHeMc2000}, we explain here each of the algorithmic components of MGRIT, starting with the coarse-grid residual equation. 

Let $\bfU$ denote an approximation to the true solution $\bfU_*$ of \eqref{eq:nonlinear_mgrit} such that $\bfU_* = \bfU + \bfE$ with $\bfE$ denoting the current error. Then this error can be expressed in terms of the residual $\bfR$ as
\begin{align}
	\bfR &:= \bfG - \bfA(\bfU,\bftheta)  = \bfA(\bfU_*, \bftheta) - \bfA(\bfU,\bftheta)\\
	     &= \bfA(\bfU + \bfE, \bftheta) - \bfA(\bfU, \bftheta). \label{eq:mgrid:residualeq}
\end{align}
In a multigrid setting, this residual equation \eqref{eq:mgrid:residualeq}
 is solved on a coarser grid such that an approximation to the error $\bfE$ can be computed more cheaply than on the fine grid. 
In linear cases, i.e. when $\bfA$ is linear in $\bfU$, the residual equation reduces to $\bfA\bfE = \bfR$ and can thus be solved for the error $\bfE$ directly. In the nonlinear case, the residual equation 
	$\bfA(\bfV,\bftheta) = \bfA(\bfU, \bftheta) + \bfR$
is solved for $\bfV$ on the coarse grid before the error can be extracted with $\bfE = \bfV - \bfU$. 

For a given time-grid discretization $t_n = nh, n=0, \dots, N$ and $h=T/N$, the coarse grid is defined by choosing a coarsening factor $c>1$ and assigning every $c$-th time point to the next coarser time-grid with $T_n = nh_{\Delta}, n=0,\dots, N_{\Delta} = N/c$, and coarse-grid spacing $h_{\Delta}=ch$.
An example of two grid levels using a coarsening factor of $c=5$ is given in Figure \ref{fig-grid}.  
The residual $\bfR$, as well as the current approximation $\bfU$ and controls $\bftheta$, are restricted to the coarse grid with injection by choosing every $c$-th time point, i.e., the restriction of $\bfU$ is 
\begin{equation}
   \bfU_\Delta = (\bfu_\Delta^0, \bfu_\Delta^1, \dots, \bfu_\Delta^{N_\Delta}), \mbox{ where } \bfu_\Delta^n = \bfu^{nc},
\end{equation}
with $\bfR_\Delta$, $\bftheta_\Delta$ defined analogously. 
Consequently, the residual equation that is to be solved on the coarse grid reads 
\begin{align}\label{eq:coarsegrid_residual}
	\bfA_{\Delta}(\bfV_{\Delta}, \bftheta_\Delta) = \bfA_{\Delta}(\bfU_{\Delta}, \bftheta_\Delta) + \bfR_{\Delta}.
\end{align}
Here, $\bfA_{\Delta}$ denotes a re-discretization of $\bfA$ on the coarse grid utilizing a coarse-grid propagator $\bfPhi_{\Delta}$, i.e., 
\begin{align}\label{eq:coarsegrid_systemmatrix}
	\bfA_{\Delta} (\bfU_{\Delta}, \bftheta_{\Delta}) :=
	\begin{pmatrix} 
		\bfu^0_{\Delta} \\
		\bfu^1_{\Delta} - \bfPhi_{\Delta}(\bfu^0_{\Delta}, \bftheta^0_{\Delta}) \\
		\vdots \\
		\bfu^{N_\Delta}_{\Delta} - \bfPhi_{\Delta}(\bfu^{N_{\Delta}-1}_{\Delta}, \bftheta^{N_{\Delta}-1}_{\Delta}) \\
	\end{pmatrix} . 
\end{align}
An obvious choice for $\bfPhi_{\Delta}$ is a re-discretization of the problem on the coarse grid,
such as by using the same propagator as on the fine grid, but with a bigger time step {size} $h_{\Delta} = ch$, thus skipping the fine-grid time points and updating only the coarse-grid points. 
For instance,
$\bfPhi$ could be a forward or backward Euler discretization with time-step size $h$, and
$\bfPhi_\Delta$ could be a forward or backward Euler discretization with time-step size $h_{\Delta} = c
h$.\footnote{{In general, the time-grid hierarchy and the corresponding coarse-grid operator $\Phi_{\Delta}$ should be chosen such that stability of the time-stepping method on each coarse time grid is ensured. In this work, the chosen hierarchy encountered
no stability issues on coarse time grids for the chosen time-stepping method. In fact, we
never observed stability issues, even on very coarse time grids.  However, a more
thorough treatment of stability for ResNets and MGRIT is beneficial and is the topic of future research.}}
In the case of forward Euler (i.e., ResNet architecture), the coarse-grid propagator $\bfPhi_\Delta$ is given by
\begin{align}
	\bfPhi_{\Delta}(\bfu_{\Delta}^n, \bftheta^n_{\Delta})	= \bfu^{n}_{\Delta} + h_{\Delta} F(\bfu_{\Delta}^n, \bftheta^n_{\Delta}).
\end{align}
 
On the {coarse} grid, the residual equation \eqref{eq:coarsegrid_residual} is solved exactly with forward substitution. Afterwards, the error approximation on the coarse grid is extracted with $\bfE_{\Delta} = \bfV_{\Delta} - \bfU_{\Delta}$. This coarse-grid error approximation is then used to correct the fine-grid approximation $\bfU$ at coarse-grid points with $\bfU^{nc} \leftarrow \bfU^{nc} + \bfE^n_{\Delta}$.

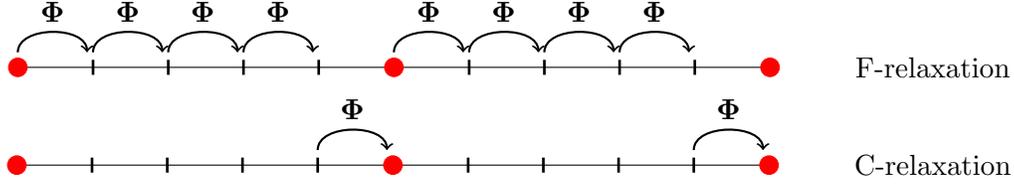
\begin{figure}
\center
\usetikzlibrary{shapes}
\usetikzlibrary{plotmarks}

\begin{tikzpicture}
\draw (0,0)--(10,0);
\filldraw (11,0) circle (0pt) node[right] {F-relaxation};
\filldraw[color = red] (0,0) circle (3.5pt);
\filldraw[color = red] (5,0) circle (3.5pt);
\filldraw[color = red] (10,0) circle (3.5pt);

\foreach \x in {1,2,3,4,6,7,8,9}
  \path[draw, line width = 1pt] (\x,0.1) -- (\x,-0.1) node {};
\foreach \x in {0,1,2,3,5,6,7,8}
 \draw[->, thick] (\x,0.2) to[out=90,in=90] node [above, midway] {$\bfPhi$} (\x+0.92,0.2); 
\end{tikzpicture}

\begin{tikzpicture}
\draw (0,0)--(10,0);
\filldraw (11,0) circle (0pt) node[right] {C-relaxation};
\filldraw[color = red] (0,0) circle (3.5pt);
\filldraw[color = red] (5,0) circle (3.5pt);
\filldraw[color = red] (10,0) circle (3.5pt);

\foreach \x in {1,2,3,4,6,7,8,9}
  \path[draw, line width = 1pt] (\x,0.1) -- (\x,-0.1) node {};
\foreach \x in {4,9}
  \draw[->, thick] (\x,0.2) to[out=90,in=90] node [above, midway] {$\bfPhi$} (\x+0.92,0.2); 
\end{tikzpicture}
  \caption{F-relaxation and C-relaxation for a coarsening by factor of $c=5$. \label{fig:FCRelax}}
\end{figure}

Complementing the coarse time-grid error correction is the fine-grid relaxation process. 
Here, block Jacobi relaxation alternates between the fine-grid and the coarse-grid points.
More precisely, relaxation on the fine points (called F-relaxation) corresponds to updating each fine point concurrently over each time chunk interval, thus propagating each coarse-point value $\bfu^{kc}$ through the corresponding fine point interval $(T^k, T^{k+1})$ as in 
\begin{align}
   \label{eq:Frelax}
   \bfu^{n} &\leftarrow \bfPhi(\bfu^{n-1},\bftheta^{n-1}), \; \mbox{ for each } n=k c +1, kc + 2, \dots, (k+1)c-1.
   \intertext{Importantly, each $k$-th interval of fine points can be computed independently, in
parallel.  Relaxation on the coarse points (called C-relaxation) is analogous, and updates 
   each coarse point concurrently by propagating the nearest left neighboring value. For the $k$-th coarse point, the update is given by}
   \label{eq:Crelax}
   \bfu^{kc} &\leftarrow \bfPhi(\bfu^{kc-1},\bftheta^{kc-1}), \; \mbox{ for } k = 1,2,\dots,N_\Delta . 
\end{align}
The actions of F- and C- relaxation are {described} in Figure~\ref{fig:FCRelax}.
Unless otherwise noted, we use FCF-relaxation, which is an
application of F-relaxation (\ref{eq:Frelax}), followed {by} an application of
C-relaxation (\ref{eq:Crelax}), and then F-relaxation (\ref{eq:Frelax}) again.
We note that such F/C orderings in relaxation are common for multigrid methods.

Taken together, the coarse-grid error correction and the fine-grid
relaxation form the two-grid MGRIT cycle depicted in Algorithm~\ref{alg:mgrit}.\footnote{The F-relaxation two-grid version of nonlinear MGRIT is equivalent to the Parareal algorithm \cite{GaVa2007}.}
Typically, the MGRIT Algorithm~\ref{alg:mgrit} is carried out recursively, with successively coarser time-grids, until a
coarsest time-grid of trivial size is reached, and Step 3 is solved exactly using forward substitution.  If the
levels are traversed in order, going down to the coarsest time-grid and then
back to the finest time-grid, this is called a V-cycle.  It corresponds to the
``Solve" in Step 3 being implemented as a single recursive call.  However, more
powerful cycles can be applied that visit coarse time-grids more frequently (such as F-cycles, see, e.g., \cite{oosterlee_book,BrHeMc2000} for more information on multigrid cycling).

Note that the main work carried out on a given time-grid is
   the parallel relaxation process.  Thus
   the work on each MGRIT level is highly parallel.  Only when a coarsest time-grid of 
   trivial size is reached, is the level solved sequentially by forward substitution.  Thus, the algorithm 
   simultaneously computes all time steps in parallel, reducing the serial propagation component to the size of the coarsest grid plus the traversal through each level.

\begin{algorithm}
   \caption{MGRIT($\bfA$, $\bfU$, $\bftheta$, $\bfG$) for two grid levels} 
\label{alg:mgrit}
\begin{algorithmic}[1]
   \STATE Apply F- or FCF-relaxation to $\bfA(\bfU, \bftheta) = \bfG$  \hfill $\triangleright$ eq. \eqref{eq:Frelax}--\eqref{eq:Crelax}
   \STATE Restrict the fine-grid approximation and residual $\bfR$ to the coarse grid:\newline
     $
     \bfU_\Delta^i \leftarrow \bfU^{ic}, \bfR_\Delta^i \leftarrow (\bfG - \bfA(\bfU, \bftheta) )^{ic},\, \text{for} \, i=0,\dots,N_\Delta
     $
	\STATE Solve $\bfA_\Delta(\bfV_\Delta, \bftheta_{\Delta}) = \bfA_\Delta(\bfU_\Delta, \bftheta_\Delta) + \bfR_\Delta$. 
			\hfill $\triangleright$ eq. \eqref{eq:coarsegrid_residual}
	\STATE Compute the coarse-grid error approximation: $\bfE_\Delta = \bfV_\Delta - \bfU_\Delta$.
	\STATE Correct $\bfU$: \newline $\bfU^{ic} \leftarrow \bfU^{ic} + \bfE_{\Delta}^{ic}$, for $i=0,\dots, N_\Delta$, and apply F-relaxation %\hfill $\triangleright$ eq. \eqref{eq:Frelax}
   \STATE \textbf{If} $\| \bfR \| \le tol$: halt. \newline
   	    \textbf{Else}: go to step 1.
\end{algorithmic}
\end{algorithm}

The MGRIT iterations can mathematically be considered a fixed-point method for solving the forward problem \eqref{eq:discretestate1}--\eqref{eq:discretestate2}. Using the iteration index $m$, it reads 
\begin{equation}
   \mbox{for } m = 0,1, \ldots: \quad \bfU_{m+1} = {\rm MGRIT}(\bfA,\bfU_m, \bftheta, \bfG).
\end{equation}
The MGRIT iterator has been shown to be a contraction in many settings for
linear, nonlinear, parabolic, and hyperbolic problems, although hyperbolic
problems tend to be more difficult (e.g., \cite{dobrev2016twolevel, Fa2014,
gunther2018xbraid, Fa2016}).
Upon convergence, the limit fixed-point $\bfU = {\rm MGRIT}(\bfA,\bfU, \bftheta,\bfG)$ will satisfy the discrete network state equations
as in \eqref{eq:discretestate1}--\eqref{eq:discretestate2}, since MGRIT solves the same underlying problem. 

Before starting the multigrid iterations, an initial solution guess for
$\bfU$ must be set.
Typically, the coarse-grid points are initialized using the best current solution
estimate.  This is often either some generic initial condition, or an
interpolated solution from a cheaper coarser time-grid.

\subsection{Multigrid Across Layers for Backpropagation} 
\label{sec:mgrit_backprop}
The same nonlinear multigrid scheme as described in Section \ref{sec:mgrit} can also be utilized to solve the adjoint equations \eqref{eq:discreteadjoint1} -- \eqref{eq:discreteadjoint2} in layer-parallel. 
The adjoint equations are linear in the adjoint variables $\bar \bfu^n$, and those are propagated backwards through the network. The adjoint space-time system thus reads

\begin{equation}
  \label{eq:adjointMGRIT}
%   A_{\bfU}(\bar \bfU, \bftheta) :=
  \underbrace{\begin{pmatrix}
		\bfI                             &     \\
		-(\partial_{\bfu}\bfPhi^{N-1})^T  & \bfI  &  \\
		   &  -(\partial_{\bfu}\bfPhi^{N-2})^T & \bfI   \\
         &    & \ddots & \ddots & \\
			&	  & &   -(\partial_{\bfu}\bfPhi^{0})^T    &    \bfI 
	\end{pmatrix}}_{=: \bfA_{\bfU}(\bar \bfU, \bftheta)}
	  \begin{pmatrix}
		\bar \bfu^{N} \\
		\bar \bfu^{N-1} \\
		\vdots \\
		\bar \bfu^1 \\
		  \bar \bfu^0
	  \end{pmatrix}
  =
	  \underbrace{\begin{pmatrix}
		\frac 1s (\partial_{\bfu^N}\ell^N)^T  \\
		0\\
		\vdots \\
		0 \\
	  \end{pmatrix}}_{ =: \bfG_{\bfU}},
\end{equation}
where again $\bar \bfu^n$ denotes the adjoint variable at layer $n$ for a general example $\bfy$ or for a batch of examples $\bfy_{k}, \, k \in\mathcal{S}\subset \{1,\dots,s\}$. Further, $(\partial_{\bfu}\bfPhi^n)^T$ denotes the partial derivative $\partial_{\bfu}\bfPhi(\bfu^n, \bftheta^n)^T {=\left(\frac{\partial \bfu^{n+1}}{\partial \bfu^n}\right)^T}$. It corresponds to the backwards layer-propagation of adjoint sensitivities which in the case of a forward Euler discretization for $\bfPhi$ (i.e., ResNet architecture), reads 
\begin{align}
	\partial_{\bfu}\bfPhi(\bfu^n, \bftheta^n)^T\bar \bfu^{n+1} = \bar \bfu^{n+1} + h\partial_{\bfu}F(\bfu^n, \bftheta^n)^T\bar \bfu^{n+1}.
\end{align}
Each backward propagator at layer $n$ depends on the primal state $\bfu^{n}$, hence the system matrix and right-hand-side of \eqref{eq:adjointMGRIT} depend on the current state $\bfU$ which is reflected in the subscript $A_{\bfU}$ and $\bfG_{\bfU}$. 
The structure of the adjoint system \eqref{eq:adjointMGRIT}, however, is the same as that of the state system~\eqref{eq:nonlinear_mgrit}. Hence the same MGRIT approach as presented in Algorithm \ref{alg:mgrit} can be utilized to solve the adjoint equations with the layer-parallel multigrid scheme by applying the following iteration
		 \begin{align}
			\bar \bfU_{m+1} = {\rm MGRIT}(\bfA_{\bfU}, \bar \bfU_m, \bftheta, \bfG_{\bfU});
		 \end{align}
for the adjoint vector $\bar{\bfU} := (\bar \bfu^N,\dots, \bar \bfu^0)$.

\begin{remark}
The adjoint equations depend on the primal states $\bfu^{n}_k$. Therefore, those states need to be either stored during forward propagation, or recomputed while solving the adjoint equations. Hybrid approaches like the check-pointing method have been developed, which compromise memory consumption with computational complexity (see, e.g., \cite{wang2009checkpointing}). Memory-free methods using reversible networks were first proposed for general dynamics in~\cite{GomezEtAl2017}. However, as shown in \cite{Chang2017Reversible}, not all architectures that are reversible algebraically are forward and backward stable numerically. This motivates limiting the forward propagation to stable dynamics, e.g., inspired by hyperbolic systems.
\end{remark}

\subsection{Non-intrusive implementation}
\label{sec:nonintrusive}

The MGRIT algorithm relies on the action of the layer-to-layer forward and backward propagators, $\bfPhi$ and $\partial_{\bfu}\bfPhi^T$, and their respective rediscretizations, $\bfPhi_\Delta$ and $\partial_{\bfu}\bfPhi_{\Delta}^T$, on coarser grid levels. However, it does not access or ``know" the internals of these functions. Hence, MGRIT can be applied in a fully non-intrusive way with respect to any existing discretization of the nonlinear dynamics describing the network forward and backward propagation. 
A user can wrap existing sequential
evolution operators according to an MGRIT software interface, and then the MGRIT
code iteratively computes the solution to \eqref{eq:discretestate1}--\eqref{eq:discretestate2} and \eqref{eq:discreteadjoint1}--\eqref{eq:discreteadjoint2} in parallel. 

Our chosen MGRIT implementation for time-parallel computations
(forwards and backwards) is XBraid \cite{xbraid-package}. 
 One particular advantage of XBraid is its generic and flexible user-interface that requires relatively straightforward user-routines which likely already exist, such
as how to take inner-products and norms with vectors $\bfu^n$, how to take
a time step with $\bfPhi$ and $\partial_{\bfu}\bfPhi^T$, etc.  

Since the user defines the action of $\bfPhi$, any existing
implementation of layer computations can continue to be used, including accelerator code, e.g., for GPUs. 
However since $\bfPhi$ takes a single time step, any use of GPU kernels for 
$\bfPhi$ implies memory movement to and from the CPU every time step.  This is
because current architectures largely rely on the CPU to handle the message
passing layer of parallelism, and it is over this layer that XBraid provides temporal
parallelism. However, future implementations could move the message
passing layer to occur solely on the GPU, thus removing this memory movement
overhead.  Additionally, the bandwidth and latency between CPUs and accelerators 
will continue to improve, also ameliorating this issue.

\begin{remark}The state and adjoint MGRIT iterations recover at convergence the same reduced gradient as a layer-serial forward- and backpropagation through the network. They can thus be integrated into any gradient-based training algorithm for updating the network control parameters $\bftheta, \bfW, \bfmu$. Sub-gradient methods, such as SGD or other batch approaches, can also be utilized by choosing the corresponding subset $\mathcal{S}\subset\{1,\dots,s\}$. 
Regarding speedup and parallelism, the layer-parallel computations are particularly attractive in the small-batch mode when options for data parallelism are limited. Overall, we expect a runtime speedup over a layer-serial approach for deep networks through the greater concurrency within the state and adjoint solves, when the computational resources are large enough. 
\end{remark}

%%%%%%%%%%%%%%%%%%%%%%%%%%%%%%%%%%%%%%%%%%%%%%%%%%%%%%%%%%%%%%%%%%%%%%
\section{Simultaneous Layer-Parallel Training}
\label{sec:oneshot}
%%%%%%%%%%%%%%%%%%%%%%%%%%%%%%%%%%%%%%%%%%%%%%%%%%%%%%%%%%%%%%%%%%%%%%
The iterative nature of the layer-parallel multigrid scheme allows for a simultaneous training approach that solves the network state and adjoint equations inexactly during training. 
To this end, we reduce the accuracy of the state and adjoint MGRIT solver during training and update the network control parameters utilizing inexact gradient information. This corresponds to an early stopping of the MGRIT iterations in each outer optimization cycle.
The theoretical background of this early-stopping approach of the inner state and adjoint fixed-point iterations is based on the \textit{One-shot} method \cite{bosse2014oneshot}, which has been successful for reducing runtimes of many PDE-constrained optimization problems in aerodynamics applications (e.g. \cite{ito2010approximate, gauger2009singlestep, bosse2014optimal}).

\begin{algorithm}
  \caption{Simultaneous Layer-Parallel Training}\label{alg:oneshot}
  \begin{algorithmic}[1]
	\STATE Perform $m_1$ state updates: \hfill  $\triangleright$ Sec. \ref{sec:mgrit}
	%  \newline
	%  \begin{flalign*}
	 	$\text{\textbf{for }} m=1,\dots, m_1: \quad \bfU_m \leftarrow {\rm MGRIT}(\bfA,\bfU_{m-1}, \bftheta,\bfG) $
	%  \end{flalign*}

	\STATE Perform $m_2$ adjoint updates: \hfill $\triangleright$ Sec. \ref{sec:mgrit_backprop}
	%  \begin{flalign*}
	 $\text{\textbf{for }} m=1,\dots, m_2: \quad \bar \bfU_{m} \leftarrow {\rm MGRIT}(\bfA_{\bfU_{m_1}},\bar \bfU_{m-1}, \bftheta, \bfG_{\bfU_{m_1}}) $
	%  \end{flalign*}

	\STATE Assemble reduced gradient $\nabla_{\bftheta}J, \nabla_{\bfW}J, \nabla_{\bfmu}J$ \hfill $\triangleright$ \eqref{eq:reducedgradient1},\eqref{eq:reducedgradient2}
			\STATE Approximate Hessians $\bfB_{\bftheta}, \bfB_{\bfW}, \bfB_{\bfmu}$ and select a stepsize $\alpha>0$ 
	\STATE Network control parameter update: \label{alg:oneshot:controlupdate}
		\newline
		$\bftheta \hspace{1.3ex} \leftarrow \bftheta \hspace{1.5ex} - \alpha \bfB_{\bftheta}^{-1}\nabla_{\bftheta}J$ \\
		$\bfW                    \leftarrow \bfW     - \alpha \bfB_{\bfW}^{-1}    \nabla_{\bfW}J$   \\
		$\bfmu    \hspace{1.2ex} \leftarrow \bfmu   \hspace{1.2ex}  - \alpha \bfB_{\bfmu}^{-1}   \nabla_{\bfmu}J$

	\STATE \textbf{If} converged: {halt} \label{oneshot:stoppingcrit} \newline 
			\textbf{Else:} go to step 1.
	\end{algorithmic}
\end{algorithm}

The simultaneous layer-parallel training approach is summarized in Algorithm \ref{alg:oneshot}. To clarify the details of the method, the following points need to be considered:
\begin{itemize}
	\item \textit{Number of state and adjoint updates $m_1, m_2$}: 
		For ``large'' $m_1, m_2$, the algorithm recovers the same gradient, and hence the same scheme as a conventional layer-serial gradient-based training approach - however with the addition of enabled layer-parallelism, providing runtime reductions through greater concurrency.  
	Considering smaller numbers of inner MGRIT iterations, e.g. $m_1,m_2\in\{1,2\}$, further reduces the runtime of each iteration and yields the simultaneous optimization approach. In that case, 
	control parameter updates in Step \ref{alg:oneshot:controlupdate} are based on inexact gradient information utilizing the most recent state and adjoint variables $(\bfu^n_{m_1})$ and $(\bar \bfu^n_{m_2})$.

	For the extreme case $m_1=m_2=1$ {(and appropriate Hessian approximation of quasi-Newton type, see below)}, the resulting optimization iteration can mathematically be interpreted as an approximate, reduced sequential quadratic programming (rSQP) method with convergence analysis presented in \cite{ito2010approximate}. In \cite{Bosse2014Adaptive}, theoretical considerations on the choice of $m_1, m_2$ are presented, which rely on the state and adjoint residuals by searching for descent on an augmented Lagrangian function. In practice, choosing $m_1,m_2$ to be as small as $2$ has proven successful in our experience. 
	
	\item \textit{Hessian approximation}:
	In order to prove convergence of the simultaneous One-shot method on a theoretical level, the preconditioners $\bfB_{\bftheta}, \bfB_{\bfW}, \bfB_{\bfmu}$ 
	should approximate the Hessian of an augmented Lagrangian function that involves the residual of the state and adjoint equations (see \cite{bosse2014oneshot} and references therein). 
	Numerically, we approximate the Hessian through consecutive limited-memory BFGS updates based on the current reduced gradient (thus assuming that the residual term is small).
	Alternatively, one might try to approximate the Hessian with a scaled identity matrix, which drastically reduces computational complexity and has already proven successful in various applications of the One-shot method. 
	 It should be noted, that the Hessian with respect to $\bfW, \bfmu$ can be computed directly as it involves only the second derivative of the loss function $\ell$ in \eqref{eq:loss} and the regularizer $R$.

	\item \textit{Stepsize selection}: The stepsize $\alpha$ is selected through a standard line-search procedure based on the current value of the objective function, e.g. a backtracking line-search satisfying the (strong) Wolfe-condition (see, e.g., \cite{nocedal2006numerical}). 

	\item \textit{Stopping criterion}: 
	Since the One-shot method targets optimality and feasibility of the state, adjoint and control variables simultaneously, the stopping criterion should involve not only the norm of the reduced gradient, but also the norm of the state and adjoint residuals. In the context of network training, however, solving the optimization problem to high accuracy is typically not desired in order to prevent overfitting. 
	We therefore compute a validation accuracy in each iteration of the above algorithm by applying the current network controls to a separate validation data set. We terminate the training, if the current network controls produce a high validation accuracy, rather than focusing on the current residuals of the state, adjoint and gradient norms. 

\end{itemize}

%%%%%%%%%%%%%%%%%%%%%%%%%%%%%%%%%%%%%%%%%%%%%%%%%%%%%%%%%%%%%%%%%%%%%%
\section{Numerical Results}
\label{sec:numerics}
%%%%%%%%%%%%%%%%%%%%%%%%%%%%%%%%%%%%%%%%%%%%%%%%%%%%%%%%%%%%%%%%%%%%%%

We investigate the computational benefits of the simultaneous layer-parallel training approach on three test cases. 
For all test cases, our focus is on the ability to achieve speedup in training runtimes for very deep neural networks by introducing parallelism between the layers. It is likely, though not explored here, that greater combined speedups are possible by additionally using data-parallelism or parallelizing inside of each layer. Further studies are required to better understand the trade-off of distributing parallel work between layer-parallel and data-parallel.

\subsection{Test Cases}

\begin{enumerate}
\item Level set classification (\textit{Peaks example}):

As a first step, we consider the test problem suggested in~\cite{HaberRuthotto2017} for classifying grid points into five  (non-convex) level sets of a smooth nonlinear function $f\colon [-3,3]^2\to\R$ (Figure \ref{fig:testcases:peaks}).
The training data set consists of $s=5000$ randomly chosen points $\bfy_k \in[-3,3]^2, \, k=1,\ldots, s,$ and {standard basis vectors} $c_k\in\R^5$ which represent the probability that a point $\bfy_k$ belongs to level set $i\in\{1,2,3,4,5\}$.
The goal is to train a network that predicts the correct level sets for new points in $[-3,3]^2$ (validation points). 

We choose a ResNet architecture with {smoothed ReLU activation defined as}
\begin{equation}
	{\sigma(x) = \begin{cases}
		\max\{x,0\}, & |x|>0.1\\
		2\frac{1}{2} x^2 + \frac{1}{2} x + \frac{1}{40} & |x| \leq 0.1 
	\end{cases}.}
\end{equation}
 Also,  we define the linear operations $\bfK(\cdot)$ at each layer to be a dense matrix representation of the weights $\bftheta^n$. We choose a network depth of $T=5$ discretized with up to $N=2048$ layers and a network width of $8$ such that $\bfu^n\in\R^8, \forall \,n=0, \dots, N$. {In order to map the data set to the network width, we choose $L_{\rm in}$ to be a dense $\R^{8\times 2}$ matrix whoses entries are learned alongside the network parameters, followed be an initial application of the activation function.}

\item Hyperspectral image segmentation (\textit{Indian Pines}):

In this test case, we consider a soil segmentation problem based on a hyperspectral image data set. The input data consists of hypersectral bands over a single landscape in Indiana, US, ({Indian Pines data set}~\cite{indianpinesdata}) with $145\times 145$ pixels. For each pixel, the data set contains $220$ spectral reflectance bands which represent different portions of the electromagnetic spectrum in the wavelength range $0.4 - 2.5 \cdot 10^{-6}$. 
The goal is to train a network that assigns each pixel of the scene to one of $16$ class labels that represent the type of land-cover present at that pixel (such as alfalfa, corn, soybean, wheat, etc.), see Figure \ref{fig:testcases:IP}.

We use the spectral bands of $s=1000$ randomly chosen pixel points, $\bfy_k\in\R^{220}, \, k=1,\dots,s$, together with their corresponding class probability vectors $\bfc_k\in\R^{16}$ (unit vectors) for training. The network architecture is a ResNet with smoothed ReLU activation (i.e. $\sigma(x) = \max\{0,x\}$, smoothed around zero) and define the linear operations $\bfK(\cdot)$ at each layer to be a dense matrix representation of the weights $\bftheta^n$. We choose a network depth of $T=20$ discretized with up to $N=2048$ layers and a network width of $220$ channels, corresponding to the 220 reflectance bands. {The initial operator $L_{\rm in}$ is chosen to be the identity.}

\item MNIST image classification (\textit{MNIST}):

As a final example, we consider the now classic MNIST~\cite{Lecun1998Gradient} test case for classification of handwritten digits encoded in a $28\times 28$ grey scale image (Figure \ref{fig:testcases:MNIST}). Our objective for this test case is to demonstrate the scalability of the layer-parallel approach over an increasing number of layers. While we obtain reasonable validation accuracy, the objective is not to develop an optimal ResNet to solve this problem. Further, we obtained the timings below with our own straightforward implementation of convolutions, to ensure compatible layer-to-layer propagators with XBraid for our initial tests.  Future work will use a fast convolution library, which will provide a substantial speedup to both the serial and layer-parallel codes.

{For the weak scaling runs below,} we use a ResNet architecture with $\tanh$ activation and define internal layers by the linear operator $\bfK(\cdot)$ using $8$ convolution kernels of width $3${; we used similar architectures in~\cite{HaberHolthamRuthotto2017,HaberRuthotto2017}}. This yields a weight tensor at each layer of size $\mathbb{R}^{3\times3\times 8\times 8}$. The parameters to be trained are {in} $\mathbb{R}^{28\times 28}$ at each layer. {The strong scaling training tests below used $4$ convolutional kernels to reduce memory requirements.} 
The network is defined to have a depth of $T=5$ and is discretized with up to $N = 2048$ layers. {The initial operator $L_{\rm in}$ is chosen to be the identity copied over the $8$ (or $4$) convolutional kernels.}

\begin{figure}
	\center
	\begin{subfigure}{0.3\textwidth}
		\center
		\includegraphics[height=3.4cm]{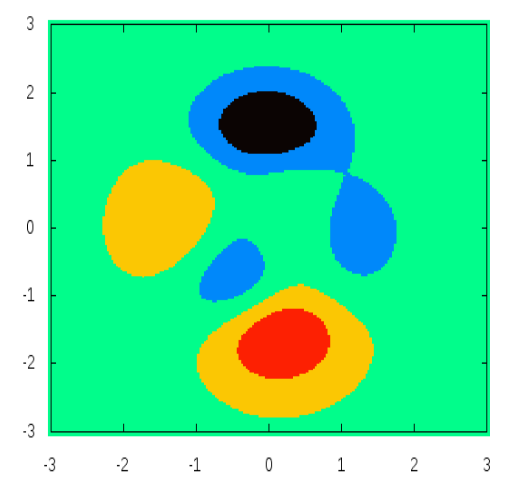}
		\caption{Peaks}
		\label{fig:testcases:peaks}
	\end{subfigure}
	\begin{subfigure}{0.3\textwidth}
	\center
		\includegraphics[height=3.4cm]{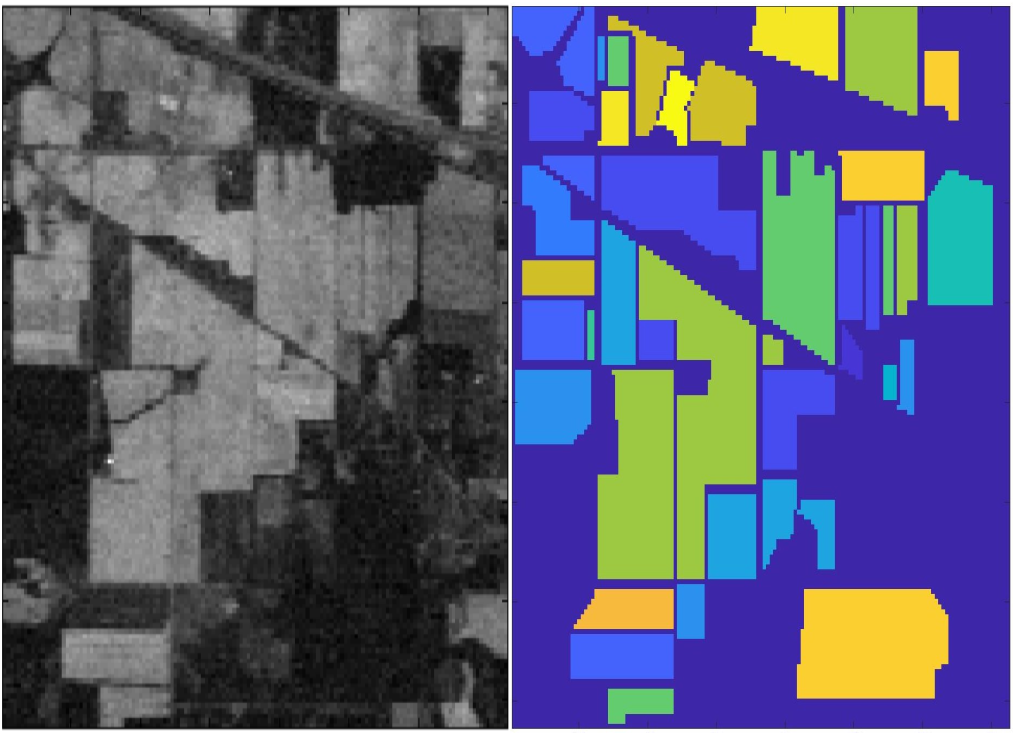} 
		\caption{Indian Pines}
		\label{fig:testcases:IP}
	\end{subfigure}
	\begin{subfigure}{0.3\textwidth}
	\center
		\includegraphics[height=3.4cm]{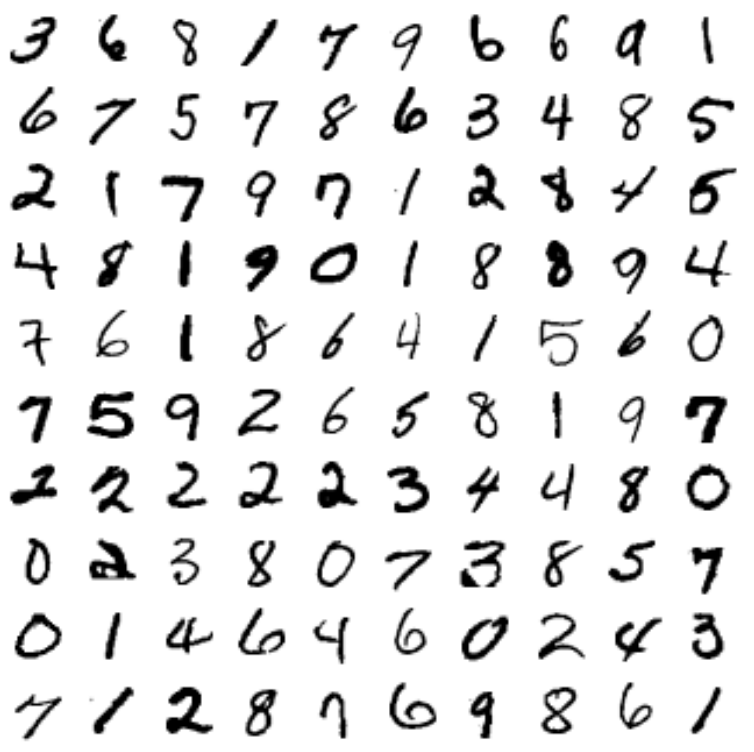}
		\caption{MNIST}
		\label{fig:testcases:MNIST}
	\end{subfigure}
	\caption{Classes of the Peaks example (test case 1), sample band and true classes of the Indian Pines data set (test case 2), and examples from the MNIST data set (test case 3).}
\end{figure}
\end{enumerate} 

The Peaks and Indian Pines computations were performed on the RHRK cluster Elwetritsch II at TU Kaiserslautern. Elwetritsch II has 485 nodes based on Haswell (2x8 cores, 64GB) and Skylake (2x12 cores, 96GB) architectures.  
The computations for the MNIST results were performed on the Skybridge capacity cluster at Sandia National Laboratories. Skybridge is a Cray containing 1848 nodes with two 8 core Intel 2.6 GHz Sandy Bridge processors, 64GB of RAM per node and an Infiniband interconnect. 
{The source code is available online at \cite{dnn_pint}}.

\subsection{Layer-Parallel Scaling and Performance Validation}
\label{sec:numerics:MGRIT}
First, we investigate {the performance of the layer-parallel MGRIT propagation for one single objective function and gradient evaluation.} Here, we keep the network weights fixed and propagate a batch of examples of sizes $s=5000, 1000, 500$ for the Peaks, Indian Pines and MNIST test case, respectively, through the network.
We choose a coarsening factor of $c=4$ to set up a hierarchy of ever coarser layer-grids to employ the multigrid scheme. {This coarsening strategy did not encounter any stability issues for forward Euler on the coarser layer-grids.} 

\begin{figure}
	\center 
	\begin{subfigure}{0.49\textwidth}
		\center 
		\includegraphics[width=\textwidth]{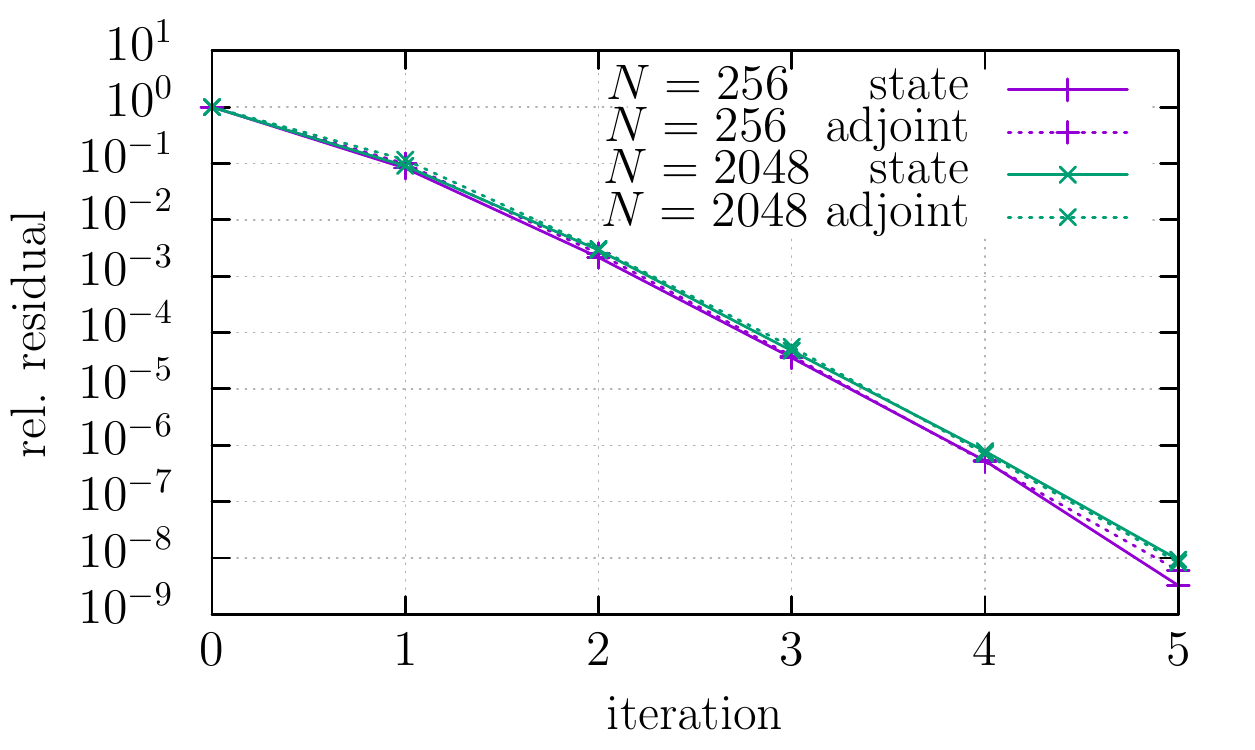}
		\caption{Peaks example}
	\end{subfigure}
	% \begin{subfigure}{0.49\textwidth}
		% \center 
		% \includegraphics[width=\textwidth]{figures/hyperspectral/convergence_standalone.pdf}
		% \caption{Indian Pines}
	% \end{subfigure}
		\begin{subfigure}{0.49\textwidth}
		\center 
		\includegraphics[width=\textwidth]{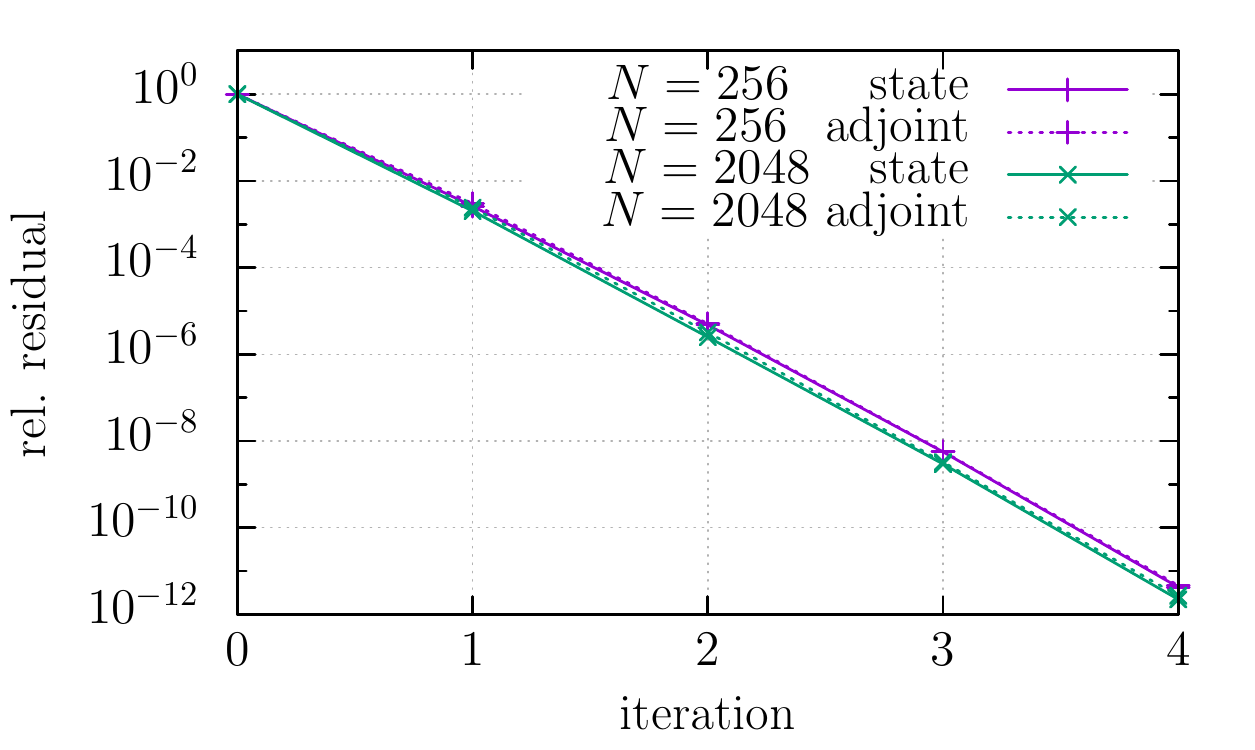}
		\caption{MNIST}
	\end{subfigure}
   \caption{Convergence history of MGRIT solving the state and adjoint equations for $N=256$ and $N=2048$ layers. The MGRIT scheme achieves fast convergence independent of the number of layers.\protect\footnotemark}
	\label{fig:MGRITconvergence}
\end{figure}

Figure \ref{fig:MGRITconvergence} shows the convergence history of the MGRIT iterations for two different problem sizes using $N=256$ and $N=2048$ layers. We monitor the relative drop of the state and adjoint residual norms and observe fast convergence for all test cases that is independent of the number of layers.
{Note that the performed multigrid iterations themselves are not dependent on the number of cores used for parallelisation, making Figure \ref{fig:MGRITconvergence} independent of the parallel distribution. We report scaling results varying the number of cores next.}
\footnotetext{The corresponding figure for the Indian Pines test case shows the same quantitative behavior, and has hence been omitted here.}

\begin{figure}
	\center
	% \begin{subfigure}{0.49\textwidth}
		% \includegraphics[width=\textwidth]{figures/peaks/xbraid/weakscaling_standalone.pdf}
		% \caption{Peaks example}
	% \end{subfigure}
	\begin{subfigure}{0.49\textwidth}
	   \includegraphics[width=\textwidth]{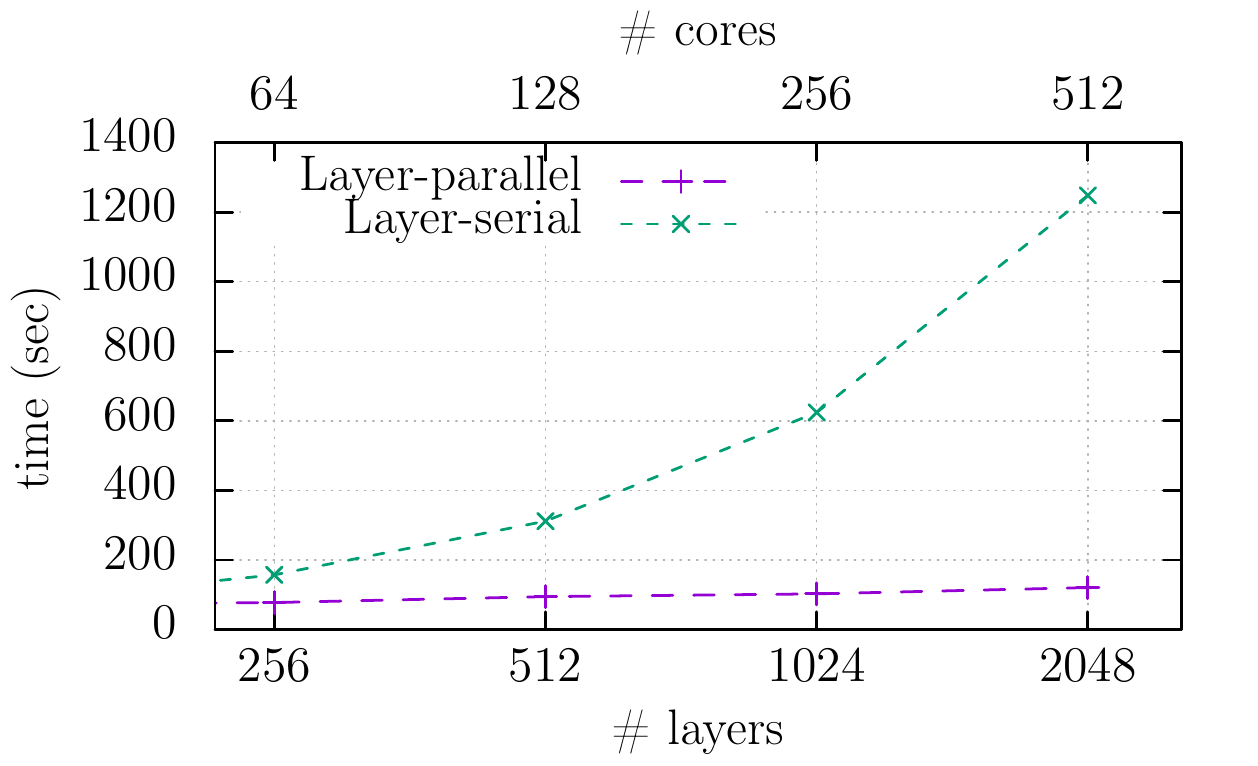}
	   \caption{Indian Pines}
	\end{subfigure}
 	\begin{subfigure}{0.49\textwidth}
	   \includegraphics[width=\textwidth]{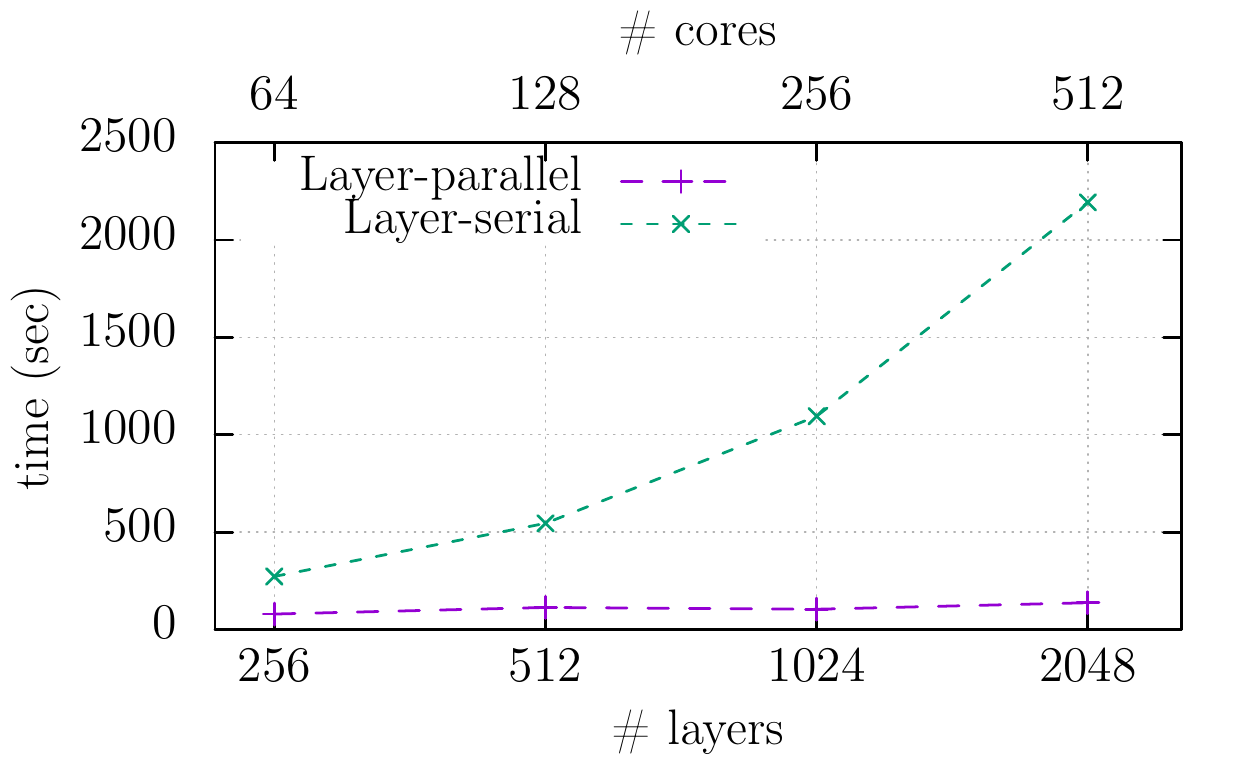}
	   \caption{MNIST}
	\end{subfigure}
	\caption{Runtime comparison of a layer-parallel gradient evaluation with layer-serial forward- and backpropagation. The layer-parallel approach yields nearly constant runtimes for increasing problem sizes and computational resources.\protect\footnotemark}
	\label{fig:weakscaling}
\end{figure}

{We investigate scaling results for the layer-parallel MGRIT scheme and compare runtimes to conventional serial-in-layer forward- and backpropagation.}
Figure \ref{fig:weakscaling} presents a weak-scaling study for the layer-parallel MGRIT scheme. Here, we double the number of layers as well as the number of compute cores while keeping the ratio $N/\#\text{cores} = 4$ fixed, such that each compute unit processes $4$ layers. Runtimes are measured for one objective function and gradient evaluation, using a relative stopping criterion of $5$ orders of magnitude for the MGRIT residual norms. Note, that the layer-serial data points have been added for comparison, even though they are executed on only one core. For the layer-serial propagation, doubling the number of layers leads to a doubling in runtime. 
The layer-parallel MGRIT approach however yields nearly constant runtimes independent of the problem size. 
The resulting speedups are reported in Table \ref{tab:MGRITspeedup}. Since the layer-parallel MGRIT approach removes the linear runtime scale of the conventional serial-layer propagation, resulting speedups increase linearly with the problem size yielding up to a factor of $16$x for the MNIST case using $2048$ layers and $512$ cores. Further speedup can be expected when considering ever more layers (and computational resources).

\footnotetext{The corresponding figure for the Peaks test case shows the same quantitative behavior, and has hence been omitted here.}

% \begin{table}
%   \center
%   \begin{tabular}{@ { } llrrrr @ { }}
%     \toprule
%      $\#$Layers    & $\#$Cores  & Peaks & Indian Pines & MNIST \\
% 	 \midrule
% 	256  & 64  &   1.8   & 2.0   &  3.4  \\
%     512  & 128 &   2.6   & 3.3   &  4.8  \\
% 	1024 & 256 &   4.9   & 6.1   &  10.5 \\
% 	2048 & 512 &   8.1   & 10.3  &  16.0 \\
%     \bottomrule
%   \end{tabular}
%   \caption{Runtime speedup of layer-parallel gradient evaluation over layer-serial propagation.}
%   \label{tab:MGRITspeedup}
% \end{table}

{
\begin{table}
  \center
  \begin{tabular}{@ { } lllrrr @ { }}
    \toprule
     Test case  & $\#$Layers    & $\#$Cores  & Serial & Parallel & Speedup\\
    %  $\#$Layers    & $\#$Cores  & serial & parallel & speedup & serial & parallel & speedup & serial & parallel & speedup \\
	 \midrule
      Peaks     & 	256  & 64  &  1.8sec  & 1.2sec   &   1.5   \\
                & 	512  & 128 &  3.7sec  & 1.5sec   &   2.5   \\
                & 	1024 & 256 &  7.1sec  & 1.6sec   &   4.3   \\
                & 	2048 & 512 & 13.9sec  & 1.8sec   &   7.7   \\
	 \midrule
	  Indian Pines & 	256  & 64  &  157.1sec  & 77.6sec   &   2.0   \\
                   & 	512  & 128 &  311.6sec  & 94.5sec   &   3.3   \\
                   & 	1024 & 256 &  624.0sec  & 102.6sec  &   6.1   \\
                   & 	2048 & 512 & 1248.0sec  & 120.6sec  &   10.3  \\
	 \midrule
	  MNIST        & 	256  & 64  &  272.3sec   &   79.5sec    &   3.4   \\
                   & 	512  & 128 &  545.3sec   &  113.3sec    &   4.8   \\
                   & 	1024 & 256 & 1095.2sec  &  104.0sec    &   10.5  \\
                   & 	2048 & 512 & 2193.5sec  &  137.3sec    &   16.0  \\	 
    \bottomrule
  \end{tabular}
  \caption{Runtime and speedup of layer-parallel gradient evaluation over layer-serial propagation.}
  \label{tab:MGRITspeedup}
\end{table}
}

\begin{figure}
	\center 
	\begin{subfigure}{.49\textwidth}
		\center 
		\includegraphics[width=\textwidth]{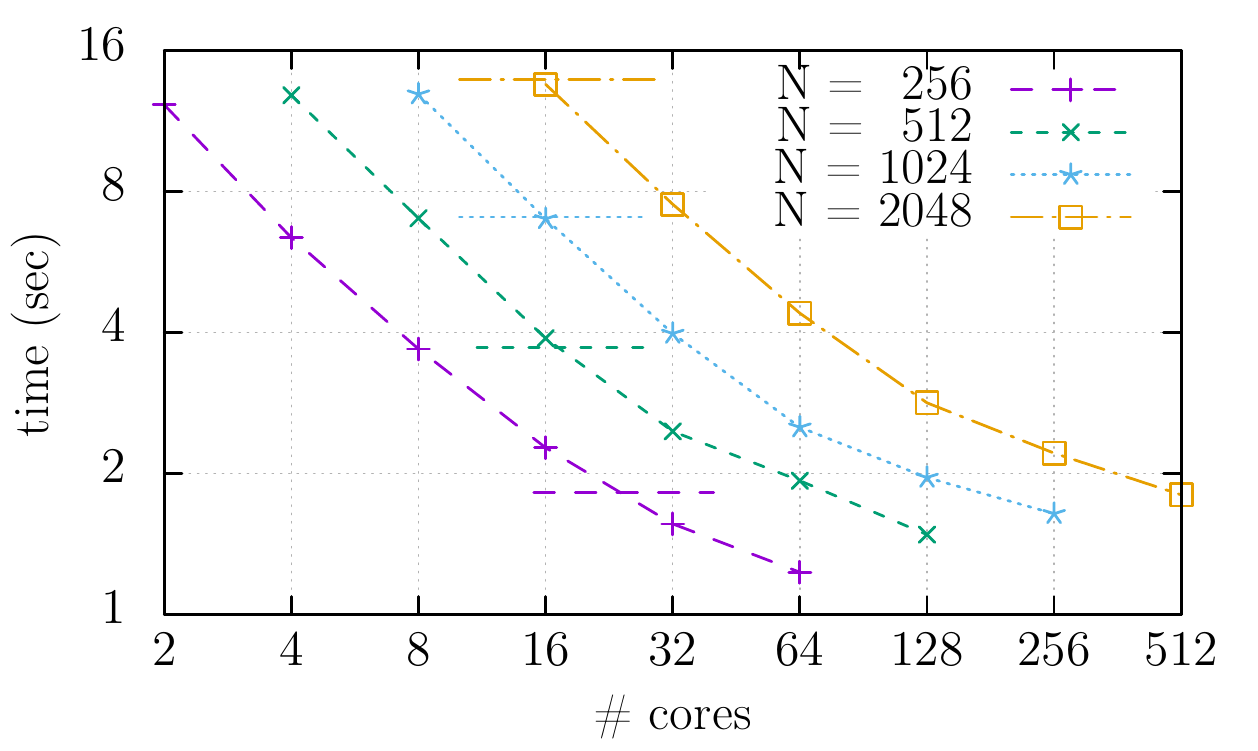}
		\caption{Peaks example}
	\end{subfigure}
	\begin{subfigure}{.49\textwidth}
		\center 
		\includegraphics[width=\textwidth]{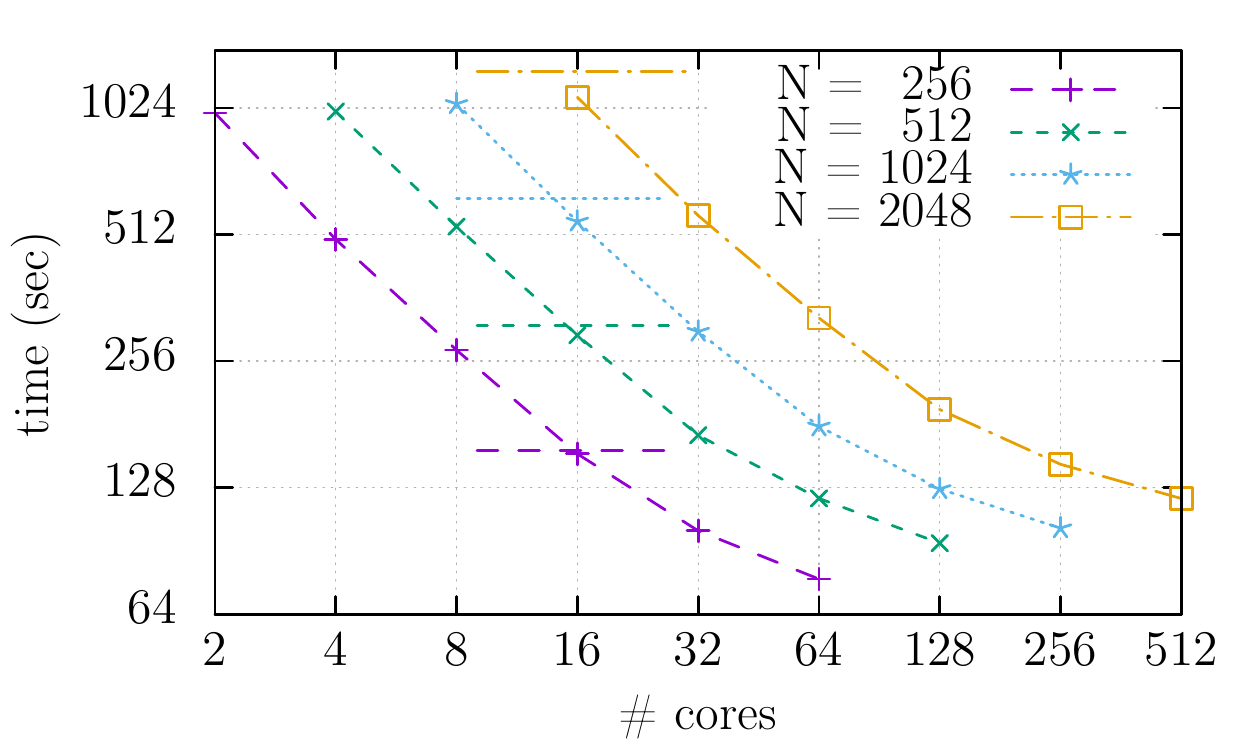}
		\caption{Indian Pines}
	\end{subfigure}
	\caption{Strong scaling study for a layer-parallel gradient evaluation for various problem sizes from $N=256$ to $N=2048$ layers. {Corresponding serial runtimes are indicated by horizontal dashed lines.} The cross-over point where the layer-parallel approach yields speedup over the layer-serial propagation lies around $16$ cores.\protect\footnotemark}
	\label{fig:strongscaling}
\end{figure}

A strong scaling study is presented Figure \ref{fig:strongscaling} for various numbers of layers. Here, we keep the problem sizes fixed and measure the time-to-solution for one gradient evaluation with MGRIT for increasing numbers of computational resources. It shows good strong scaling behavior for all test cases, independent of the numbers of layers. The cross over point where the layer-parallel MGRIT approach shows speedup over the layer-serial propagation is around $16$ cores for all cases.  

\footnotetext{The corresponding figure for the MNIST test case shows the same quantitative behavior, and has hence been omitted here.}

\subsection{Simultaneous Layer-Parallel Training Validation}
\label{sec:numerics:oneshot}
Next, we investigate the simultaneous layer-parallel training, using $m_1 = m_2 = 2$ layer-parallel MGRIT iterations in each outer training iteration (see Algorithm \ref{alg:oneshot}). {The Hessian approximations $B_{\bftheta}, B_{\bfW}, B_{\bfmu}$ are computed by successive limited-memory BFGS updates based on the current gradient $\nabla_{(\bftheta, \bfW, \bfmu)}J$.}
We compare runtimes of the simultaneous layer-parallel training with a conventional layer-serial training approach, while choosing the same Hessian, as well as the same initial network parameters for both approaches. However, we tune the optimization hyper-parameters (such as regularization parameters, stepsize selection, etc.)  separately for both schemes, in order to find the best setting for either approach that reaches a prescribed validation accuracy with the least iterations and minimum runtime. 

\begin{figure}[ht]
	\center
	\begin{subfigure}{0.49\textwidth}
		\center
		\includegraphics[width=\textwidth]{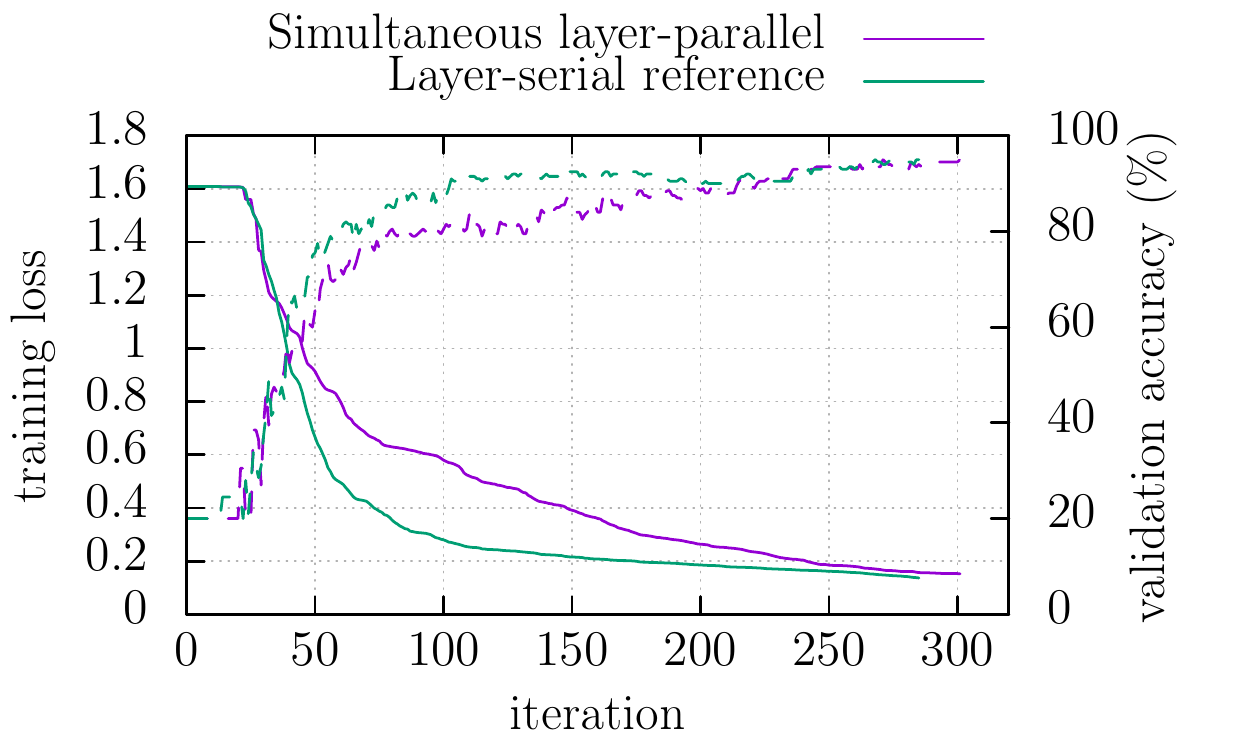}
		\caption{Peaks: Training over iteration counts}
	\end{subfigure}
	\begin{subfigure}{0.49\textwidth}
		\center
		\includegraphics[width=\textwidth]{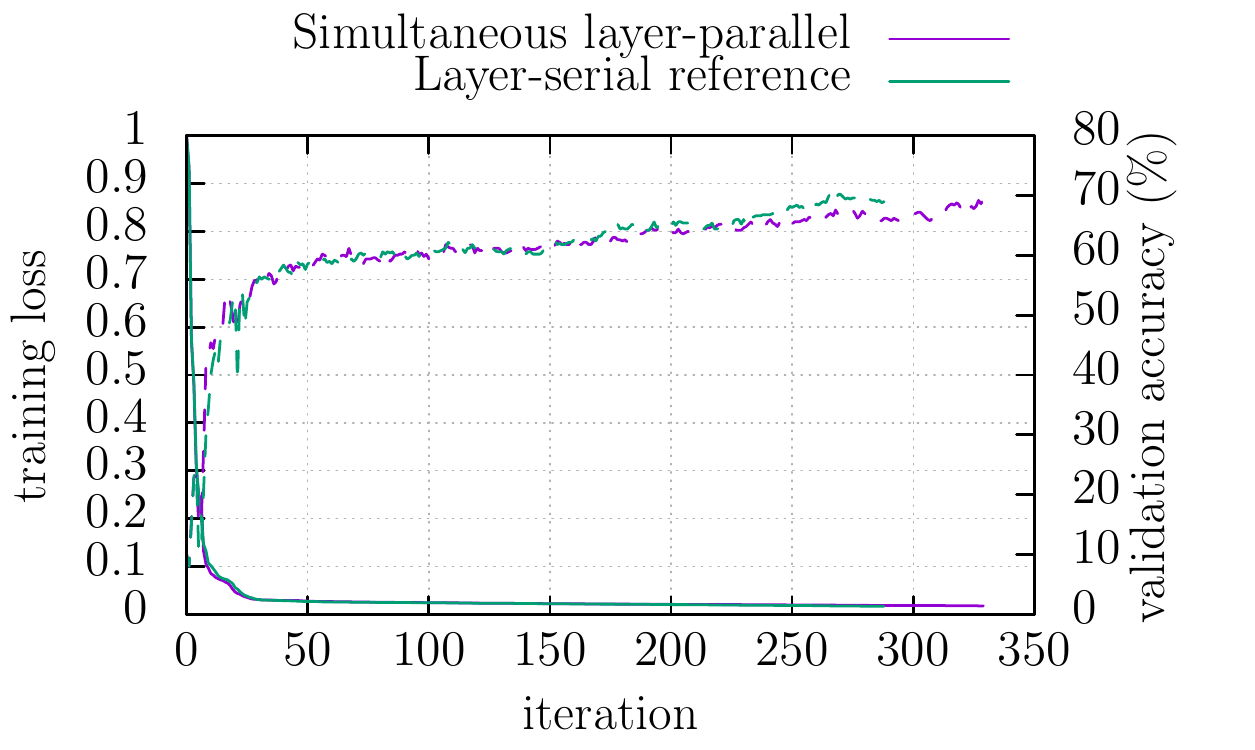}
		\caption{Indian Pines: Training over iteration counts}
	\end{subfigure}

	\begin{subfigure}{0.49\textwidth}
		\center
		\includegraphics[width=\textwidth]{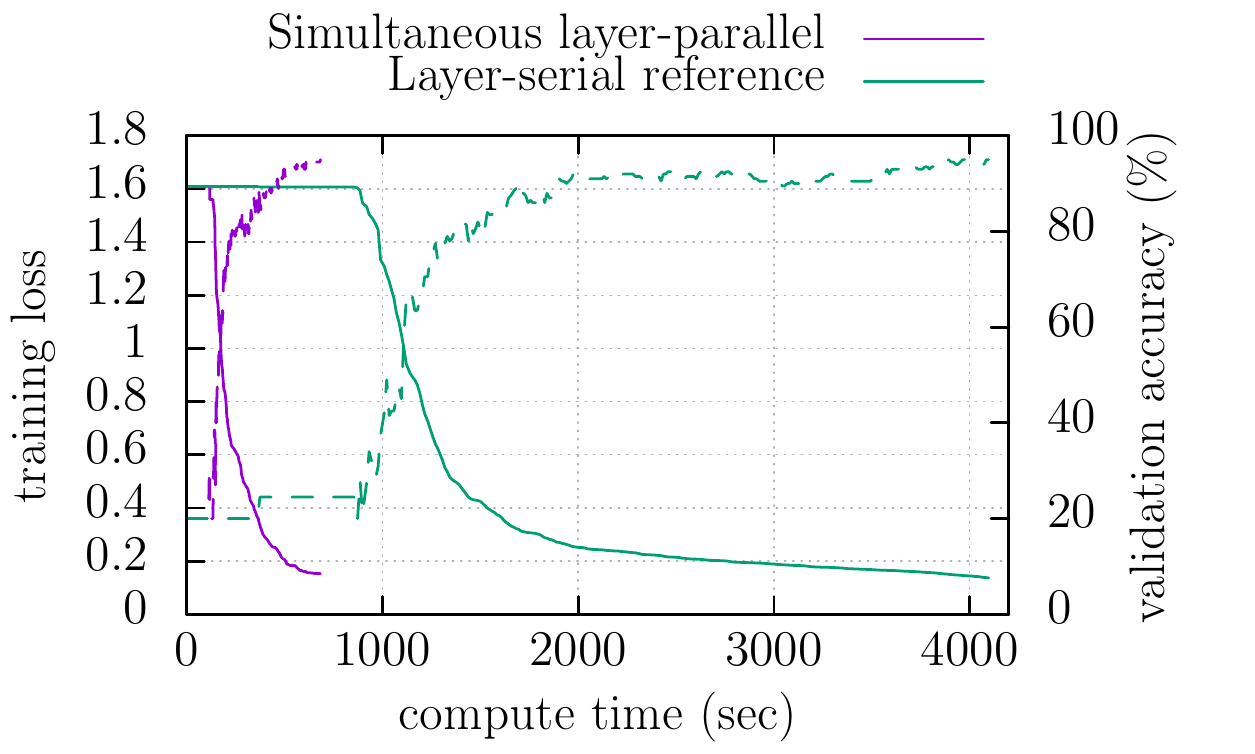}
		\caption{Peaks: Training over time}
	\end{subfigure}
	\begin{subfigure}{0.49\textwidth}
		\center
		\includegraphics[width=\textwidth]{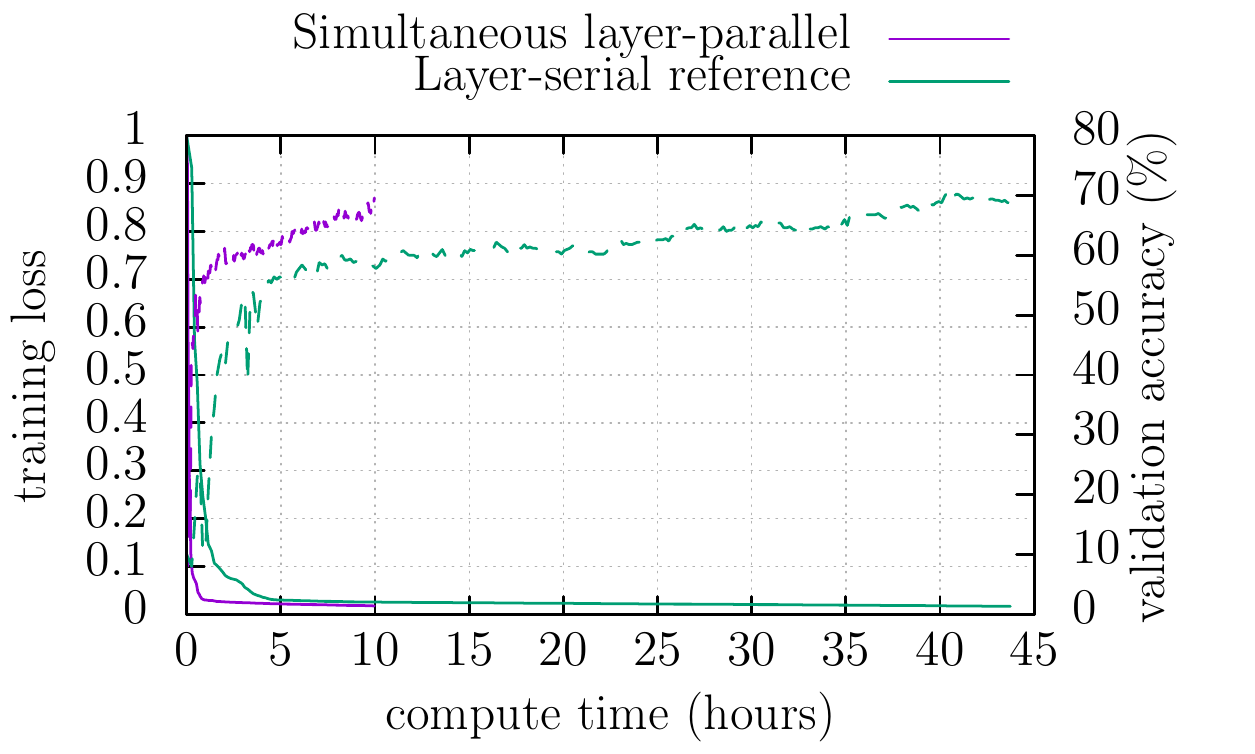}
		\caption{Indian Pines: Training over time}
	\end{subfigure}
	\caption{Training loss (solid lines) and validation accuracy (dashed lines) over training iterations (top) and compute time (bottom). For the layer-parallel training, each core processes $4$ layers. The simultaneous layer-parallel approach reaches training results comparable to a layer-serial approach within much less computational time.\protect\footnotemark}
	\label{fig:oneshot}
\end{figure}

For the Peaks example, we train a network with $N=1024$ layers distributed onto $256$ compute cores, and for the Indian Pines data set and the MNIST case we choose $N=512$ layers distributed onto $128$ compute cores, giving $4$ layers per processor in all cases. 
 Figure \ref{fig:oneshot} plots the training history over iteration counts (top) as well as runtime (bottom). We validate from the top figures, that both approaches reach comparable performance in terms of training result (optimization iteration counts, training loss and validation accuracy). Hence, reducing the accuracy of the inner multigrid iterations for solving the state and adjoint equations within a simultaneous training framework does not deteriorate the training behavior.
 But, each iteration of the simultaneous layer-parallel approach is much faster than for the layer-serial approach due to the layer-parallelization and the reduced state and adjoint accuracy. Therefore, the overall runtime for reaching that same final training result is reduced drastically (bottom figures). 
 Runtime speedups are reported in Table \ref{tab:OSspeedup}.
%  We observe a speedup of 6.0 for the Peaks example when using 256 cores, which is a parallel efficiency of about 2.3\%.
%  In the larger Indian Pines example and the MNIST case, the reported speedup of $4.4$ and $8.5$ using $128$ cores give a parallel efficiency of $3.4\%$ and $6.6\%$.
  While these results have been computed for selected fixed $N$, it is expected that the speedup scales linearly with increasing numbers of layers, similar to the observation in Table~\ref{tab:MGRITspeedup}.

\footnotetext{The corresponding figures for the MNIST test case show the same quantitative behavior, and have hence been omitted here.}

\begin{table}[h]
  \center
  \begin{tabular}{@ { } lrrrrr @ { }}
    \toprule
    Test case      & $N$  & $\#$Cores  & Layer-serial & Layer-parallel & Speedup \\
	 \midrule
	 Peaks example  & 1024 & 256       &  4096 sec      &   683 sec     &  6.0   \\
	 Indian Pines   & 512  & 128       &  2623 min      &   597 min     &   4.4  \\
	 MNIST          & 512  & 128       &  619  min       &   71  min     &   8.5  \\

    \bottomrule
  \end{tabular}
  \caption{Runtime speedup of simultaneous layer-parallel training over layer-serial training.}
  \label{tab:OSspeedup}
\end{table}

%%%%%%%%%%%%%%%%%%%%%%%%%%%%%%%%%%%%%%%%%%%%%%%%%%%%%%%%%%%%%%%%%%%%%%
\section{Conclusion}
\label{sec:conclustion}
%%%%%%%%%%%%%%%%%%%%%%%%%%%%%%%%%%%%%%%%%%%%%%%%%%%%%%%%%%%%%%%%%%%%%%
In this paper, we provide a proof-of-concept for layer-parallel training of deep residual neural networks (ResNets). 
The similarity of training ResNets to the optimal control of nonlinear time-dependent differential equations motivates us to use parallel-in-time methods that have been popular in many engineering applications.
The method developed is based on nonlinear multigrid methods and introduces a new form of parallelism across layers. 

We demonstrate two options to benefit from the layer-parallel approach. First, the nonlinear multigrid reduction in time (MGRIT) method can be used to replace forward and backward propagation in existing training algorithms, including for stochastic approximation methods such as SGD. In our experiments, this leads to speedup over serial implementations when using more than 16 compute {cores}. Second, additional savings can be obtained through the simultaneous layer-parallel training, which uses only inexact forward and backward propagations.

While the reported speedups might seem small in terms of parallel efficiency, these reductions can be of significant importance when considering large overall training runtimes.   When bare training runtimes are in the order of days, any runtime reduction is appreciated, as long as computational resources are available. Further, since training a network typically involves a careful choice of hyper-parameters,
faster training runtimes will enable faster hyper-parameter optimization and thus eventually lead to better training results in general.  
Lastly, we mention that such efficiencies for multigrid-in-time are not
uncommon \cite{falgout2017multigrid}, where the nonintrusiveness of MGRIT
contributes to the seemingly low efficiency, as does the fact that we are defining the
efficiency of MGRIT with respect to an optimal serial algorithm.  If the efficiency were
defined with respect to MGRIT using 1 core, then the efficiencies would be much higher. 

Motivated by these first promising results, we will investigate the use of layer-parallel training for more challenging learning tasks, including more complex image-recognition problems. Further reducing the memory footprint of our algorithm in those applications motivates the use of reversible networks arising from hyperbolic systems~\cite{Chang2017Reversible}. 
A challenge arising here is the interplay of MGRIT and hyperbolic systems.
{ Lastly, we  note that while the current work focused on algorithmic
development, it could nonetheless benefit greatly from an integration with a more
optimized code such as TensorFlow or Chainer.  This is planned future work.}
{An interesting topic concerning the interplay of existing codes with MGRIT is to develop strategies that handle more complicated layer architectures such as pooling and connector layers. Here, we plan on leveraging previous work on adaptive spatial coarsening/refinement in MGRIT where the problem size at each time step can change.}

\bibliography{2018-ParallelInTimeDNN} 
\bibliographystyle{abbrv}

\end{document}